\documentclass[11pt]{amsart}
\usepackage{mathbbol}
\usepackage{amsmath, amsthm, amssymb}

\newtheorem{thm}{Theorem}[section]
\newcommand{\bt}{\begin{thm}}
\newcommand{\et}{\end{thm}}

\newtheorem{cor}[thm]{Corollary}
\newcommand{\bc}{\begin{cor}}
\newcommand{\ec}{\end{cor}}

\newtheorem{lem}[thm]{Lemma}
\newcommand{\bl}{\begin{lem}}
\newcommand{\el}{\end{lem}}

\newtheorem{prop}[thm]{Proposition}
\newcommand{\bp}{\begin{prop}}
\newcommand{\ep}{\end{prop}}

\newtheorem{defn}[thm]{Definition}
\newcommand{\bd}{\begin{defn}}      
\newcommand{\ed}{\end{defn}}

\newtheorem{rmrk}[thm]{Remark}
\newcommand{\br}{\begin{rmrk}}
\newcommand{\er}{\end{rmrk}}

\newcommand{\thmref}[1]{Theorem~\ref{#1}}
\newcommand{\secref}[1]{Section~\ref{#1}}
\newcommand{\lemref}[1]{Lemma~\ref{#1}}

\newcommand{\corref}[1]{Corollary~\ref{#1}}
\newcommand{\propref}[1]{Proposition~\ref{#1}}
\newcommand{\remref}[1]{Remark~\ref{#1}}

\newcommand{\N}{\mathbb{N}}

\newcommand{\R}{\mathbb{R}}
\newcommand{\Z}{\mathbb{Z}}
\newcommand{\Hyp}{\mathbb{H}}
\newcommand{\dist}{\operatorname{dist}}
\newcommand{\diam}{\operatorname{diam}}

\newcommand{\hm}{{\mathcal H}}
\newcommand{\lm}{{\mathcal L}}
\newcommand{\md}{\operatorname{md}}

\newcommand{\jac}{{\mathbf J}}
\newcommand{\spaceform}{\mathcal{M}}

\newcommand{\lip}{\operatorname{Lip}}
\newcommand{\mass}[2][]{{\mathbf M_{#1}}(#2)}

\newcommand{\form}{{\mathcal D}}        
\newcommand{\curr}{{\mathbf M}}         
\newcommand{\norcurr}{{\mathbf N}}      
\newcommand{\intrectcurr}{{\mathcal I}} 
\newcommand{\intcurr}{{\mathbf I}}      

\newcommand{\fillvol}{{\operatorname{Fillvol}}}

\newcommand{\rstr}{\:\mbox{\rule{0.1ex}{1.2ex}\rule{1.1ex}{0.1ex}}\:}
\newcommand{\bdry}{\partial\hspace{-0.05cm}}
\newcommand{\slice}[3]{\langle#1,#2,#3\rangle}
\newcommand{\on}[1]{|_{#1}}
\newcommand{\spt}{\operatorname{spt}}
\newcommand{\ohne}{\backslash}

\newcommand{\cone}[2]{{#1\:\mbox{$\times\hspace{-0.6em}\times$}\:#2}}

\newcommand{\divergence}{{\operatorname{div}}}
\newcommand{\divGMT}{{\widehat{\operatorname{div}}}}
\newcommand{\Erank}{{\operatorname{Rank}}}

\begin{document}

\title[Filling invariants at infinity]{Filling invariants at infinity and the Euclidean rank of Hadamard spaces}

\author{Stefan Wenger}

\address
  {Courant Institute of Mathematical Sciences\\
   251 Mercer Street\\
   New York, NY 10012}
\email{wenger@cims.nyu.edu}

\date{July 25, 2006}

\begin{abstract}
In this paper we study a homological version of the asymptotic filling invariant $\divergence_k$ defined by Brady and Farb in \cite{Brady-Farb} and 
show that it is a quasi-isometry invariant for all proper cocompact Hadamard spaces, i.e.~proper cocompact
${\rm CAT}(0)$-spaces, and that it can furthermore be used to detect the Euclidean rank of such spaces.
We thereby extend results of \cite{Brady-Farb, Leuzinger, Hindawi-filling-invariants} from the setting of symmetric spaces of non-compact type to that of 
Hadamard spaces.
Finally, we exhibit the optimal growth of the $k$-th homological divergence
for symmetric spaces of non-compact type with Euclidean rank no larger than $k$ and for ${\rm CAT}(\kappa)$-spaces with $\kappa<0$.
\end{abstract}

\keywords{Filling invariants, Hadamard spaces, Alexandrov spaces, Euclidean rank, integral currents}



\maketitle

\bigskip

\section{Introduction and statement of the main results}
\sloppy

In \cite{Brady-Farb} Brady and Farb introduced a new quasi-isometry invariant $\divergence_k(X)$ of a cocompact Hadamard manifold $X$. 
The $k$-dimensional divergence $\divergence_k(X)$ can be seen as a higher dimensional analogue of the divergence of geodesics 
(studied by Gersten in \cite{Gersten-graph-manifolds} and 
\cite{Gersten-CAT}) and in some sense measures the $(k+1)$-dimensional spread of geodesic rays in $X$.
Brady and Farb then proved that $\divergence_{k-1}(X)$ has exponential growth when $X=\Hyp^{m_1}\times\dots\times\Hyp^{m_k}$ is the product of $k$ hyperbolic planes.
The main idea in their proof was to construct a family of quasi-isometric embeddings of $\Hyp^{m_1+\dots+m_k-k+1}$ in $X$ transversal to a maximal flat.
Using the same idea Leuzinger \cite{Leuzinger} extended this result to symmetric spaces of non-compact type. In particular, he showed that
$\divergence_k(X)$ grows exponentially for $k=\Erank X -1$, where $\Erank X$ denotes the Euclidean rank of $X$, i.e.~the maximal dimension of an isometrically embedded 
Euclidean space in $X$. In
\cite{Brady-Farb} the authors asked whether $\divergence_k(X)$ can be used to detect the Euclidean rank $\Erank X$ of a symmetric space $X$.
This question has recently been settled in the affirmative by Hindawi \cite{Hindawi-filling-invariants} who
showed that $\divergence_k(X)$ grows polynomially of degree at most $r^{k+1}$ for $k\geq \Erank X$. 
The primary aim of this article is to define a homological version of $\divergence_k(X)$ and show that it can be used to detect the Euclidean rank of all 
proper cocompact Hadamard spaces $X$.

The homological divergence $\divGMT_k(X)$ is defined using integral currents instead of Lipschitz images of Euclidean balls as was done in \cite{Brady-Farb}.
We use the theory of  integral currents in metric spaces developed by Ambrosio and Kirchheim in \cite{Ambr-Kirch-curr}. Precise definitions will be given 
in \secref{Section:Preliminaries}.
Roughly speaking, $\divGMT_k(X)$ is a three-parameter family of functions 
$\delta^k_{x_0,\varrho, A}$ where $x_0\in X$, $A>0$, and $0<\varrho\leq 1$. For fixed parameters, $\delta^k_{x_0,\varrho, A}(r)$ is defined to be the maximal mass of 
an integral $(k+1)$-current with support outside the open ball $U(x_0,\varrho r)$ of radius $\varrho r$ around $x_0$ 
needed to fill a $k$-dimensional integral cycle with compact support in the metric sphere of radius $r$ around $x_0$ 
and with mass at most $Ar^k$. 

In our first result we show that $\divGMT_k$ is a quasi-isometry invariant for all proper cocompact Hadamard spaces. This in particular generalizes 
\cite[Theorem 1.1]{Brady-Farb} from the context of cocompact Hadamard manifolds to that of proper cocompact Hadamard spaces.

\bt\label{theorem:quasi-isometry-invariance}
 Let $X$ and $Y$ be proper, cocompact, and quasi-isometric Hadamard spaces and $k\in\N$. Then 
 \begin{equation*}
  \divGMT_k(X)\sim_{k+1} \divGMT_k(Y),
 \end{equation*} 
 i.e.~the divergence functions $\divGMT_k(X)$ and $\divGMT_k(Y)$ have the same asymptotic growth up to an additive term  of the form $cr^{k+1}$.
\et
This theorem follows from a more general result proved in \secref{section:quasi-isometry-invariance}. For the definition of $\sim_{k+1}$ see also
\secref{section:quasi-isometry-invariance}.

For our second and main result we will need the following terminology concerning different types of growth.
\bd\label{Definition:growth-type}
 Let $X$ be a complete metric space, $k\in\N$, and $\beta\in[1,\infty)$. We write $\divGMT_k(X)\preceq r^\beta$ if there exist $0<\varrho_0\leq 1$ 
 such that
 \begin{equation*}\label{equation:growth-div-delta}
  \limsup_{r\to\infty} \frac{\delta^k_{x_0,\varrho_0, A}(r)}{r^\beta}<\infty
 \end{equation*}
 for all $x_0\in X$ and all $A>0$. 
 On the other hand, we write $\divGMT_k(X)\succeq r^\beta$ if there exists $A_0>0$ such that
 \begin{equation*}
  \liminf_{r\to\infty} \frac{\delta^k_{x_0,\varrho, A_0}(r)}{r^\beta}>0
 \end{equation*}
 for all $x_0\in X$ and all $\varrho\in(0,1]$.
\ed
%

The main result of the paper can now be stated as follows.
\bt\label{Theorem:rank-div}
Let $X$ be a proper cocompact Hadamard space and let $k\in\N$. If $k=\Erank X - 1$ then $\divGMT_k(X)\succeq r^{k+2}$. On the other hand, if
$k\geq\Erank X$ then $\divGMT_k(X)\preceq r^{k+1}$.
\et
The second part of the theorem generalizes the main result (Theorem 1.1) in \cite{Hindawi-filling-invariants} from the context of symmetric
spaces of non-compact type to that of proper cocompact Hadamard spaces. The first part of the theorem is in the spirit of the results in \cite{Brady-Farb, Leuzinger}.
It should be mentioned here that our methods of proof are different from those in \cite{Brady-Farb}, \cite{Leuzinger}, and \cite{Hindawi-filling-invariants}. In our
approach we use the isoperimetric inequality proved in \cite{Wenger-GAFA}, the characterization of the Euclidean rank via asymptotic cones 
\cite{Kleiner-local-structure}, and techniques from geometric measure theory in metric spaces. 
An important tool from the latter will be the equivalence of weak and flat convergence for integral currents proved in \cite{Wenger-flatconv}.
 
A direct consequence of the above theorem is the following corollary.
\bc\label{Corollary:detect-rank}
 The $\divGMT_k$ can be used to detect the Euclidean rank of all proper cocompact Hadamard spaces.
\ec
%
%
The estimates in \thmref{Theorem:rank-div} are good enough to detect the Euclidean rank of every proper cocompact Hadamard space
but the optimal growth rate of $\divGMT_k$ is believed to be different.\\*[1ex]
{\bf Question:} Let $(X,d)$ be a proper cocompact Hadamard space. Is it true that
 \begin{enumerate}
  \item the divergence $\divGMT_k(X)$ grows exponentially if $k=\Erank X - 1$?
  \item the divergence $\divGMT_k(X)$ grows polynomial of degree $k$ if $k\geq \Erank X$?
 \end{enumerate}
As mentioned above, for symmetric spaces of non-compact type (i) follows from \cite{Leuzinger}.
In the following we give an answer to (ii) when $X$ is a symmetric space of non-compact type or a complete ${\rm CAT}(\kappa)$-space with $\kappa<0$. These results 
are consequences of a simple relation between $\divGMT_k(X)$ and the type of isoperimetric inequality in $X$.
To state the results we adopt the following terminology.
\bd
 Let $k\in\N$, $\alpha\in\left[1,\frac{k+1}{k}\right]$ and let $X$ be a complete metric space. We say that $X$ admits an isoperimetric inequality of power 
 $\alpha$ for $\intcurr_k(X)$ if there exists a constant $C$ such that for every $T\in\intcurr_k(X)$ with $\bdry T=0$ there is an $S\in\intcurr_{k+1}(X)$ with
 $\bdry S=T$ and
\begin{equation*}
 \mass{S}\leq C[\mass{T}]^\alpha.
\end{equation*}
\ed
In the above definition, $\intcurr_k(X)$ denotes the space of $k$-dimensional metric integral currents introduced in \cite{Ambr-Kirch-curr}. Furthermore, $\mass{T}$ 
is the mass of $T$ and $\bdry T$ its boundary. See \secref{Section:Preliminaries} for the definitions.
%
%

In \cite{Gromov-filling} Gromov proved that every Hadamard manifold $X$ admits an isoperimetric inequality of at least Euclidean type for $\intcurr_k(X)$, 
thus with power at most $\alpha:=\frac{k+1}{k}$. In \cite{Wenger-GAFA} Gromov's result was extended to arbitrary Hadamard spaces.
As for ${\rm CAT}(\kappa)$-spaces with $\kappa<0$ we can quite easily prove the following linear isoperimetric inequality.
\bt\label{Theorem:linear-isoperimetric-inequality}
 Let $(X,d)$ be a complete ${\rm CAT}(\kappa)$-space with $\kappa<0$. Then for every $k\geq 1$ and every $T\in\intcurr_k(X)$
 with $\bdry T=0$ and of bounded support there exists $S\in\intcurr_{k+1}(X)$ with $\bdry S= T$ and such that
 \begin{equation*}
  \mass{S}\leq \frac{1}{\sqrt{-\kappa}k}\mass{T}.
 \end{equation*}
\et
This theorem comes as a consequence of a more general result, \thmref{theorem:cat-cone-inequality}, which can also be used to derive a monotonicity formula 
for minimizing currents in Hadamard spaces, \corref{corollary:hadamard-monotonicity-formula}. 
Theorems~\ref{Theorem:linear-isoperimetric-inequality}, \ref{theorem:cat-cone-inequality} and 
\corref{corollary:hadamard-monotonicity-formula} were previously known for simply-connected Riemannian manifolds with the suitable upper curvature bound.

We point out that a conjecture of Gromov (somewhat implicitly contained in \cite{Gromov-asymptotic}) states that every proper cocompact Hadamard space admits
a linear isoperimetric inequality for $\intcurr_k(X)$ when $k\geq\Erank X$. This conjecture is open for most classes of spaces. It is known to be true for symmetric
spaces of non-compact type and it is clearly true for complete ${\rm CAT}(\kappa)$-spaces with $\kappa<0$, as shows \thmref{Theorem:linear-isoperimetric-inequality}.
An affirmative answer to Gromov's conjecture together with the following simple relation between the isoperimetric inequality and the growth of $\divGMT_k(X)$
would immediately answer Question (ii) above.
\bp\label{Proposition:isop-div}
 Let $X$ be a Hadamard space and $k\in\N$. If $X$ admits an isoperimetric inequality of power $\alpha<\frac{k+1}{k}$ for $\intcurr_k(X)$ 
 then $\divGMT_k(X)\preceq r^{\alpha k}$.
\ep 
We can now exhibit the optimal growth rate for $\divGMT_k(X)$ for symmetric spaces $X$ of non-compact type and for complete 
${\rm CAT}(\kappa)$-spaces $X$ with $\kappa<0$ and answer question (ii) for these spaces. This generalizes Theorems 1.1 and 1.4 in \cite{Hindawi-filling-invariants}.
\bc\label{Corollary:div-CAT-negative}
 If $X$ is a complete ${\rm CAT}(\kappa)$-space with $\kappa<0$ then $\divGMT_k(X)\preceq r^k$ for all $k\in\N$.
\ec
This is an immediate consequence of \propref{Proposition:isop-div} and \thmref{Theorem:linear-isoperimetric-inequality}.
\bc
 If $X$ is a symmetric space of non-compact type then $\divGMT_k(X)\preceq r^k$ for all $k\geq\Erank(X)$.
\ec
This follows from \propref{Proposition:isop-div} and the well-known fact that Gromov's conjecture holds true for symmetric spaces $X$ of non-compact type.
The latter is a consequence of the fact that the orthogonal projection onto maximal flats in $X$ decreases the 
$k$-dimensional volume exponentially with the distance to the flat. This can be used to construct linear fillings.
\medskip

The structure of the paper is as follows: In \secref{Section:Preliminaries} we recall the necessary definitions concerning Hadamard spaces and 
integral currents and define the homological divergence $\divGMT_k$. \secref{Section:proof} contains the proof of \thmref{Theorem:rank-div}. 
The main purpose of \secref{Section:div-isop} is to establish the optimal cone inequality from which follow \thmref{Theorem:linear-isoperimetric-inequality}
and the monotonicity formula for minimizing integral currents, \corref{corollary:hadamard-monotonicity-formula}. The purpose of 
\secref{section:quasi-isometry-invariance} is to prove \thmref{theorem:generalized-quasi-isometry-invariance} from which \thmref{theorem:quasi-isometry-invariance}
will follow.

\bigskip

{\bf Acknowledgments:} The author would like to thank Urs Lang and Bruce Kleiner for inspiring discussions on the topic. 
The author also thanks the Erwin Schroedinger Institute ESI in Vienna for its hospitality and its financial support during a 3-month stay in 2005 where
parts of the paper were written.

\section{Preliminaries}\label{Section:Preliminaries}
The purpose of this section is to fix notation regarding ${\rm CAT}(\kappa)$-spaces and metric integral currents on the one hand and to define the 
homological divergence $\divGMT_k(X)$ on the other hand.

\subsection{Metric spaces of bounded curvature and asymptotic cones}
For a general reference on metric spaces of curvature bounded above in the sense of Alexandrov we refer the reader to \cite{Ballmann}, \cite{Bridson-Haefliger},
and \cite{Burago-Burago-Ivanov}. The notation we use in this article is consistent with that in \cite{Bridson-Haefliger}.

Let $\kappa\in\R$ and set $D_\kappa:= \frac{\pi}{\sqrt{\kappa}}$ if $\kappa>0$ and $D_\kappa:= \infty$ otherwise. We note that $D_\kappa=\diam(\spaceform^2_\kappa)$
where $\spaceform^2_\kappa$ is the $2$-dimensional simply-connected Riemannian manifold of constant sectional curvature $\kappa$.
A metric space $(X,d)$ is called ${\rm CAT}(\kappa)$ if the following two properties hold:
\begin{enumerate}
 \item $X$ is $D_\kappa$-geodesic: Any two points $x,y\in X$ with $d(x,y)<D_\kappa$ can be joined by a geodesic, i.e.~a curve of length $d(x,y)$.
 \item Every geodesic triangle in $X$ of perimeter $<2D_\kappa$ satisfies the ${\rm CAT}(\kappa)$-inequality, i.e.~it is at least as slim as a comparison triangle 
  in $\spaceform^2_\kappa$.
\end{enumerate}
We refer the reader to \cite[Definition II.1.1]{Bridson-Haefliger} for the precise definition of (ii). Following \cite{Ballmann} we call complete 
${\rm CAT}(0)$-spaces Hadamard spaces. Hadamard manifolds, i.e.~simply-connected Riemannian manifolds of non-positive sectional curvature, are examples of Hadamard spaces.
Furthermore, a metric space $X$ is called Alexandrov space of curvature bounded from above 
by $\kappa$ if for every point $x\in X$ there is a closed ball $B(x, r)$ which is ${\rm CAT}(\kappa)$. 
In the following we will write $B(x,r)$ for the closed ball $\{x'\in X:d(x,x')\leq r\}$ and $U(x,r)$ for the open ball $\{x'\in X:d(x,x')<r\}$. Furthermore, 
$S(x,r)$ will denote the metric sphere $\{x'\in X: d(x,x')=r\}$.

Recall that $X$ is said to be cocompact if there exists a compact set $K\subset X$ such that $X=\bigcup_{g\in\Gamma}gK$ where $\Gamma$ denotes the isometry
group of $X$. Furthermore, $X$ is said to be proper if every closed ball of finite radius is compact.
\bd 
 The Euclidean rank of a cocompact Hadamard space $X$ is defined to be
 \begin{equation*}
  \Erank X:= \sup\left\{n\in\N: \text{There exists an isometric embedding $\R^n\hookrightarrow X$}\right\}
 \end{equation*}
 where $\R^n$ is endowed with the Euclidean metric.
\ed

We finally recall the notion of asymptotic cone which will be used in the proof of the main result. As a general reference we mention \cite{Kleiner-Leeb}.
A non-principal ultrafilter on $\N$ is a finitely additive probability measure $\omega$ on $\N$ together with the $\sigma$-algebra of all subsets such
that $\omega$ takes values in $\{0,1\}$ only and $\omega(A)=0$ whenever $A\subset \N$ is finite.
Using Zorn's lemma it is not difficult to establish the existence of non-principal ultrafilters on $\N$, see e.g.~Exercise I.5.48 in \cite{Bridson-Haefliger}. 
It is also easy to prove the following fact.
If $(Y,\tau)$ is a compact topological Hausdorff space then for every sequence $(y_m)_{m\in\N}\subset Y$ there exists a unique point $y\in Y$ such that
\begin{equation*}
 \omega(\{m\in\N: y_m\in U\})=1
\end{equation*}
for every $U\in\tau$ containing $y$. We will denote this point by $\lim\nolimits_\omega y_m$.

Let now $(X,d)$ be a metric space and fix a non-principal ultrafilter $\omega$ on $\N$, a basepoint $\star\in X$ and a sequence $r_m\nearrow\infty$. 
Define an equivalence relation on the set of sequences $(x_m)_{m\in\N}\subset X$ satisfying
\begin{equation}\label{equation:bounded-sequence-def}
 \sup_m\frac{1}{r_m}d(\star, x_m)<\infty
\end{equation}
by
\begin{equation*}
 (x_m)\sim (x'_m)\quad\text{if and only if}\quad \lim\nolimits_\omega \frac{1}{r_m}d(x_n,x'_n)=0.
\end{equation*}
\bd
The asymptotic cone $(X, r_m^{-1}d, \star)_\omega$ is the set of equivalence classes of sequences $(x_m)\subset X$ satisfying \eqref{equation:bounded-sequence-def}
together with the metric given by
\begin{equation*}
 d_\omega([(x_m)],[(x'_m)]):= \lim\nolimits_\omega \frac{1}{r_m}d(x_m,x'_m).
\end{equation*}
\ed

We refer the reader to \cite{Kleiner-Leeb} for properties of asymptotic cones of Hadamard spaces and to \cite{Kleiner-local-structure} for the connection to 
the Euclidean rank which will be exploit in the proof of \thmref{Theorem:rank-div}.

\subsection{Lipschitz maps into metric spaces}
Let $(X,d)$ be a metric space, $U\subset\R^k$ open, and let $\varphi: U\to X$ be a Lipschitz map. In \cite[Theorem 2]{Kirchheim} Kirchheim 
proved that for almost every $z\in U$ the metric derivative
\begin{equation*}
 \md\varphi_z(v):= \lim_{r\searrow 0}\frac{d(\varphi(z+rv),\varphi(z))}{r}
\end{equation*}
exists for every $v\in\R^k$ and
\begin{equation}\label{equation:metric-differential-strong}
 \lim_{r\searrow 0}\frac{1}{r}d(\varphi(z+rv),\varphi(z+rw)) = \md\varphi_z(v-w)
\end{equation}
for all $v,w\in\R^k$. This was independently discovered by Korevaar and Schoen, see \cite{Korevaar-Schoen}.
It follows from \eqref{equation:metric-differential-strong} that $\md\varphi_z$ is a seminorm on $\R^k$ for almost every $z\in U$ and, if $\varphi$ is bi-Lipschitz, that
it is even a norm.
If $U\subset\R^k$ is merely measurable then $\md\varphi_z$ can be defined at almost every Lebesgue density point $z\in U$ by a simple approximation argument. 

The $k$-th Jacobian of a semi-norm $s$ on $\R^k$ is defined by
\begin{equation*}
 \jac_k(s):= \frac{\omega_k}{\lm^k(\{v\in\R^k: s(v)\leq 1\})},
\end{equation*}
where $\lm^k$ denotes the Lebesgue measure on $\R^k$ and $\omega_k$ is the Lebesgue measure of the unit ball in $\R^k$.
If $s$ is a norm then
\begin{equation*}
 \jac_k(s) = \frac{\hm^k_s(Q)}{\lm^k(Q)}
\end{equation*}
whenever $Q\subset\R^k$ has strictly positive measure. 
Here, $\hm^k_s$ denotes the $k$-dimensional Hausdorff measure on $(\R^k, s)$. We recall that given a metric space $(Y,d)$ the Hausdorff measure of
$A\subset Y$ is defined by
\begin{equation*}
 \hm_Y^k(A):=\lim_{\delta\searrow 0}\inf\left\{\sum_{i=1}^\infty\omega_k \left(\frac{\diam(B_i)}{2}\right)^k: 
  A\subset\bigcup_{i=1}^\infty B_i, \diam(B_i)<\delta\right\}.
\end{equation*}

Finally, the area factor of a norm $s$ on $\R^k$ is given by
\begin{equation*}
 \lambda_s:= \max\left\{\frac{\lm^k(L([0,1]^k))}{\jac_k(s)}: \text{ $L=(L_1,\dots, L_k): (\R^k,s)\to\R^k$ lin, $L_i$ $1$-lip}\right\}.
\end{equation*}
If $k=1$ of if $s$ comes from an inner product then $\lambda_s=1$. We point out that for a normed space $(V,\|\cdot\|)$ the $k$-volume density
\begin{equation*}
 \mu(v_1\wedge\dots\wedge v_k):= \lambda_{\|\cdot\|_{\operatorname{span}\{v_1,\dots,v_k\}}}\hm^k_V(v_1\wedge\dots\wedge v_k)
\end{equation*}
is usually called Gromov mass$^*$ or Benson volume density, see e.g.~\cite{Alvarez-Thompson}.

\subsection{Integral currents in metric spaces}
The general reference for this section is the work of Ambrosio and Kirchheim \cite{Ambr-Kirch-curr} where the theory of normal and integral currents is extended
from the setting of Euclidean space to arbitrary complete metric spaces. The classical Euclidean theory was developed to a large part by Federer and Fleming, 
see \cite{Federer-Fleming} and \cite{Federer}.

Let $(X,d)$ be a complete metric space and $k\geq 0$ and let $\form^k(X)$ denote the set of $(k+1)$-tuples $(f,\pi_1,\dots,\pi_k)$ 
of Lipschitz functions on $X$ with $f$ bounded. The Lipschitz constant of a Lipschitz function $f$ on $X$ will
be denoted by $\lip(f)$.
\bd
A $k$-dimensional metric current  $T$ on $X$ is a multi-linear functional on $\form^k(X)$ satisfying the following
properties:
\begin{enumerate}
 \item If $\pi^j_i$ converges point-wise to $\pi_i$ as $j\to\infty$ and if $\sup_{i,j}\lip(\pi^j_i)<\infty$ then
       \begin{equation*}
         T(f,\pi^j_1,\dots,\pi^j_k) \longrightarrow T(f,\pi_1,\dots,\pi_k).
       \end{equation*}
 \item If $\{x\in X:f(x)\not=0\}$ is contained in the union $\bigcup_{i=1}^kB_i$ of Borel sets $B_i$ and if $\pi_i$ is constant 
       on $B_i$ then
       \begin{equation*}
         T(f,\pi_1,\dots,\pi_k)=0.
       \end{equation*}
 \item There exists a finite Borel measure $\mu$ on $X$ such that
       \begin{equation}\label{equation:mass-def}
        |T(f,\pi_1,\dots,\pi_k)|\leq \prod_{i=1}^k\lip(\pi_i)\int_X|f|d\mu
       \end{equation}
       for all $(f,\pi_1,\dots,\pi_k)\in\form^k(X)$.
\end{enumerate}
\ed
The space of $k$-dimensional metric currents on $X$ is denoted by $\curr_k(X)$ and the minimal Borel measure $\mu$
satisfying \eqref{equation:mass-def} is called mass of $T$ and written as $\|T\|$. We also call mass of $T$ the number $\|T\|(X)$ 
which we denote by $\mass{T}$.
The support of $T$ is, by definition, the closed set $\spt T$ of points $x\in X$ such that $\|T\|(B(x,r))>0$ for all $r>0$. 
\br
As is done in \cite{Ambr-Kirch-curr} we will also assume here that the cardinality of any set is an Ulam number. This is consistent with the 
standard ZFC set theory. We then have that $\spt T$ is separable and furthermore that $\|T\|$ is concentrated on a 
$\sigma$-compact set, i.\ e.\ $\|T\|(X\ohne C) = 0$ for a $\sigma$-compact set $C\subset X$ (see \cite{Ambr-Kirch-curr}).
\er
The restriction of $T\in\curr_k(X)$ to a Borel set $A\subset X$ is given by 
\begin{equation*}
  (T\rstr A)(f,\pi_1,\dots,\pi_k):= T(f\chi_A,\pi_1,\dots,\pi_k).
\end{equation*}
This expression is well-defined since $T$ can be extended to a functional on tuples for which the first argument lies in 
$L^\infty(X,\|T\|)$.\\
The boundary of $T\in\curr_k(X)$ is the functional
\begin{equation*}
 \bdry T(f,\pi_1,\dots,\pi_{k-1}):= T(1,f,\pi_1,\dots,\pi_{k-1}).
\end{equation*}
It is clear that $\bdry T$ satisfies conditions (i) and (ii) in the above definition. If $k=0$ or if $\bdry T$ also has 
finite mass (condition (iii)) then $T$ is called a normal current. The respective space is denoted by $\norcurr_k(X)$.\\
The push-forward of $T\in\curr_k(X)$ 
under a Lipschitz map $\varphi$ from $X$ to another complete metric space $Y$ is given by
\begin{equation*}
 \varphi_\# T(g,\tau_1,\dots,\tau_k):= T(g\circ\varphi, \tau_1\circ\varphi,\dots,\tau_k\circ\varphi)
\end{equation*}
for $(g,\tau_1,\dots,\tau_k)\in\form^k(Y)$. This defines a $k$-dimensional current on $Y$, as is easily verified.\\
The basic example of a $k$-dimensional metric current on $\R^k$ is given by
\begin{equation*}
 \Lbrack\theta\Rbrack(f,\pi_1,\dots,\pi_k):=\int_K\theta f\det\left(\frac{\partial\pi_i}{\partial x_j}\right)\,d\lm^k
\end{equation*}
for all $(f,\pi_1,\dots,\pi_k)\in\form^k(\R^k)$, where $K\subset\R^k$ is measurable and $\theta\in L^1(K,\R)$.\\
In this paper we will mainly be concerned with integral currents. We recall that an $\hm^k$-measurable set $A\subset X$
is said to be countably $\hm^k$-rectifiable if there exist countably many Lipschitz maps $\varphi_i :B_i\longrightarrow X$ from subsets
$B_i\subset \R^k$ such that
\begin{equation*}
\hm^k\left(A\ohne \bigcup \varphi_i(B_i)\right)=0.
\end{equation*}

An element $T\in\curr_0(X)$ is called integer rectifiable if there exist finitely many points $x_1,\dots,x_n\in X$ and $\theta_1,\dots,\theta_n\in\Z\ohne\{0\}$ such
that
\begin{equation*}
 T(f)=\sum_{i=1}^n\theta_if(x_i)
\end{equation*}
for all bounded Lipschitz functions $f$.
 A current $T\in\curr_k(X)$ with $k\geq 1$ is said to be integer rectifiable if the following properties hold:
 \begin{enumerate}
  \item $\|T\|$ is concentrated on a countably $\hm^k$-rectifiable set and vanishes on $\hm^k$-negligible Borel sets.
  \item For any Lipschitz map $\varphi:X\to\R^k$ and any open set $U\subset X$ there exists $\theta\in L^1(\R^k,\Z)$ such that 
    $\varphi_\#(T\rstr U)=\Lbrack\theta\Rbrack$.
 \end{enumerate}
The space of integer rectifiable currents is denoted by $\intrectcurr_k(X)$.
Integer rectifiable normal currents are called integral currents. The corresponding space is denoted by $\intcurr_k(X)$.
As the mass of a $k$-dimensional normal current vanishes on $\hm^k$-negligible sets (\cite[Theorem 3.7]{Ambr-Kirch-curr}) it 
is easily verified that the push-forward of an integral current under a Lipschitz map is again an integral current. 
It was shown in \cite[Theorem 11.1]{Ambr-Kirch-curr} that in case $X=\R^n$ the spaces $\norcurr_k(X)$, respectively $\intrectcurr_k(X)$ and $\intcurr_k(X)$, are in 
one-to-one correspondence with the spaces of normal, respectively integer rectifiable and integral currents defined by Federer and Fleming. Furthermore, 
integer rectifiable currents in a complete metric space $X$ can be represented as countable sums of 
$\varphi_{i\#}\Lbrack\theta_i\Rbrack$ where $\theta_i\in L^1(K_i,\Z)$ with $K_i\subset\R^k$ compact, $\varphi_i:K_i\to X$ bi-Lipschitz and $\varphi_i(K_i)$ pairwise 
disjoint as was shown in \cite[Theorem 4.5]{Ambr-Kirch-curr}. Moreover,
\begin{equation*}
 \mass{\varphi_{i\#}\Lbrack\theta_i\Rbrack}=
  \int_{K_i}|\theta_i(z)|\lambda_{\md\varphi_{iz}}\jac_k(\md\varphi_{iz})d\lm^k(z).
\end{equation*}

In the following, an element $T\in\intcurr_k(X)$ with zero boundary $\bdry T=0$ will be called a cycle. An element $S\in\intcurr_{k+1}(X)$
satisfying $\bdry S=T$ is said to be a filling of $T$.

We end this section with the following product construction defined in \cite{Wenger-GAFA}. It is a straight-forward generalization 
of the cone construction given in \cite{Ambr-Kirch-curr}.
For this endow $[0,1]\times X$ with the Euclidean product metric and let $f\in\lip([0,1]\times X)$. For $x\in X$ and $t\in[0,1]$ we write $f_t(x):=f(t,x)$.
With $T\in\norcurr_k(X)$ and $t\in[0,1]$ we associate a $k$-dimensional normal current on $[0,1]\times X$ by
\begin{equation*}
 ([t]\times T)(f,\pi_1,\dots,\pi_k):=T(f_t,\pi_{1\,t},\dots,\pi_{k\,t})
\end{equation*}
for $(f,\pi_1,\dots,\pi_k)\in\form^k([0,1]\times X)$. We also associate with $T$ the functional
\begin{equation*}
 \begin{split}
 ([0,1]&\times T)(f,\pi_1,\dots,\pi_{k+1})\\ &:=\sum_{i=1}^{k+1}(-1)^{i+1}\int_0^1T\left(f_t\frac{\partial\pi_{i\,t}}{\partial t},\pi_{1\,t},\dots,\pi_{i-1\,t},\pi_{i+1\,t},
  \dots,\pi_{k+1\,t}\right)\,dt
\end{split}
\end{equation*}
for $(f,\pi_1,\dots,\pi_{k+1})\in\form^{k+1}([0,1]\times X)$. It can be checked (see \cite{Ambr-Kirch-curr} and \cite{Wenger-GAFA}) that 
$[0,1]\times T\in\norcurr_{k+1}([0,1]\times X)$ and
\begin{equation}\label{equation:product-boundary}
 \bdry([0,1]\times T)=[1]\times T - [0]\times T - [0,1]\times\bdry T.
\end{equation}
If $T\in\intcurr_k(X)$ then furthermore $[0,1]\times T\in\intcurr_{k+1}(X)$.

\subsection{The $k$-th homological divergence ${\mathbf \divGMT_k(X)}$}
 As above, let $(X,d)$ be a complete metric space and $k\geq 0$.
\bd
 The $k$-th homological divergence $\divGMT_k(X)$ of $X$ is the three parameter family
 \begin{equation*}
  \divGMT_k(X):= \left\{ \delta^k_{x_0, \varrho, A}: x_0\in X, 0<\varrho\leq 1, A>0\right\}
 \end{equation*}
 where $\delta^k_{x_0, \varrho, A}(r)$ is the function given by
 \begin{equation*}
  \begin{split}
   \delta^k_{x_0, \varrho, A}(r)= \sup\Big\{\fillvol_{X\ohne U(x_0, \varrho r)}(T):\; T\in\intcurr_k(X),\, &\spt T\subset S(x_0,r)\text{ {\rm cpt}},\\
           &\bdry T=0,\, \mass{T}\leq Ar^k\Big\}
  \end{split}
 \end{equation*}
 if $k\geq 1$ and
 \begin{equation*}
  \begin{split}
   \delta^0_{x_0, \varrho, A}(r)= \sup\Big\{\fillvol_{X\ohne U(x_0, \varrho r)}(T):\; T\in\intcurr_0(X),\, &\spt T\subset S(x_0,r),\\
           &T(1)=0,\, \mass{T}\leq A\Big\}.
  \end{split}
 \end{equation*}
\ed
Here, for a closed subset $C\subset X$ and a $T\in\intcurr_k(X)$ with $\spt T\subset C$ the filling volume of $T$ in $C$ is defined as
\begin{equation*}
 \fillvol_C(T):= \inf\left\{\mass{S}: S\in\intcurr_{k+1}(X), \bdry S=T, \spt S\subset C\right\}
\end{equation*}
where we agree on $\inf\emptyset=\infty$. We note that if $X$ is such that any two points can be joined by a Lipschitz curve then for $T\in\intcurr_0(X)$ the 
condition $T(1)=0$ is equivalent to the condition that there exists an $S\in\intcurr_1(X)$ with $\bdry S=T$.

Given $f,g:[0,\infty)\to[0,\infty)$ and $\beta\in[1,\infty)$ we write $f\preceq_\beta g$ if there exist constants $a,b,c>0$ such that $f(r)\leq ag(br)+cr^\beta$ for 
all $r>0$ large enough.
We furthermore write $f\sim_\beta g$ if $f\preceq_\beta g$ and $g\preceq_\beta f$. This defines an equivalence relation on functions from $[0,\infty)$ to $[0,\infty)$.
\bl\label{Lemma:basepoint-delta}
 Let $(X,d)$ be a Hadamard space and $A>0$, $\varrho\in(0,1)$. Then $\delta^k_{x_0, \varrho, A}\preceq_k\delta^k_{x_0', \varrho', A'}$ for all $x_0, x'_0\in X$, 
 all $\varrho'>\varrho$ and $A'>A$.
\el
\begin{proof}
 Fix $\varrho, \varrho', A, A'$ as in the hypothesis. For $x_0, x'_0\in X$ set $L:=d(x_0,x'_0)$. 
 Let $T\in\intcurr_k(X)$ with $\bdry T=0$ and with $\spt T\subset S(x_0, r)$ compact and $\mass{T}\leq Ar^k$, where
 $r$ is large enough. Set $r':= r- L$ and denote by $\pi:X\to B(x'_0,r')$ the orthogonal projection, see \cite[Proposition II.2.4]{Bridson-Haefliger}. 
 Define $\varphi: [0,1]\times X\to X$ to be
 the locally Lipschitz map for which $t\mapsto \varphi(t,x)$ is the constant speed geodesic from $\pi(x)$ to $x$. Set $S':= \varphi_{\#}([0,1]\times T)$ and note 
 that $S'\in\intcurr_{k+1}(X)$ and, by \eqref{equation:product-boundary}, we furthermore have $\bdry S' = T- \pi_{\#}T$. For 
 $(f,\tau_1,\dots,\tau_{k+1})\in\form^{k+1}(X)$ we abbreviate $$(\hat{f}_t,\hat{\tau}_{1\,t},\dots,\hat{\tau}_{k+1\, t}):= 
 (f\circ\varphi_t,\tau_1\circ\varphi_t,\dots,\tau_{k+1}\circ\varphi_t)$$ and compute
 \begin{align*}
  \big|S'(f,\tau_1,\dots,&\tau_{k+1})\big|\\ 
  &\leq \sum_{i=1}^{k+1}\left|\int_0^1T\Big(\hat{f}_t\,\frac{\partial\hat{\tau}_{i\,t}}{\partial t}, 
   \hat{\tau}_{1\,t},\dots,\hat{\tau}_{i-1\,t},\hat{\tau}_{i+1\, t}, \dots,\hat{\tau}_{k+1\,t}\Big)dt\right|\\
  &\leq \sum_{i=1}^{k+1}\int_0^1\prod_{j\not= i}\lip(\hat{\tau}_{j\,t})\int_X\Big|\hat{f}_t\,\frac{\partial\hat{\tau}_{i\,t}}{\partial t}\Big|\,d\|T\|\,dt\\
  &\leq (k+1)L\prod_{j= 1}^{k+1}\lip(\tau_j)\int_0^1\int_X|f\circ\varphi(x,t)|\,d\|T\|\,dt.
 \end{align*}
From this it follows that $\|S'\|\leq (k+1)L\varphi_{\#}(\lm^1\times\|T\|)$ and, in particular, $$\mass{S'}\leq (k+1)L\mass{T}\leq (k+1)LAr^k.$$
If $r$ is chosen large enough we furthermore have 
\begin{equation*}
 \mass{\pi_\#T}\leq \mass{T}\leq Ar^k\leq A'r'^{\,k}.
\end{equation*}
Since $\pi_\# T$ has compact support in $S(x'_0,r')$ there exists by definition an $S''\in\intcurr_{k+1}(X)$ with compact support 
satisfying $\bdry S''=\pi_\# T$ as well as 
$\mass{S''}\leq \delta^k_{x'_0,\varrho',A'}(r)$ and $\spt S''\subset X\ohne U(x'_0,\varrho'r')$.
It follows that $S:= S' + S''\in\intcurr_{k+1}(X)$ satisfies $\bdry S=T$ and 
\begin{equation*}
 \mass{S}\leq \delta^k_{x'_0,\varrho',A'}(r)+L(k+1)Ar^k.
\end{equation*}
Furthermore $\spt S$ is compact and, if $r$ is chosen large enough, $\spt S\subset X\ohne U(x_0, \varrho r)$. This completes the proof.
\end{proof}

\section{Proof of the main result}\label{Section:proof}
The aim of this section is to prove \thmref{Theorem:rank-div}. We will need the following lemma which is a variation of 
an argument of Ambrosio and Kirchheim \cite[Theorem 10.6]{Ambr-Kirch-curr}, 
see also \cite[Lemma 3.4]{Wenger-GAFA}. It yields the existence of fillings with sufficient volume growth.
\bl\label{Lemma:Q-min-growth}
 Let $X$ be a Hadamard space, $k\in\N$, and $\alpha\in\left[1,\frac{k+1}{k}\right]$. Suppose $X$ admits an isoperimetric inequality
 of power $\alpha$ for $\intcurr_k(X)$ and set $C:=\max\{C',C''\}$ where $C'$ and $C''$ are the constants of the isoperimetric inequalities of power $\alpha$
and $\frac{k+1}{k}$, respectively. Then there exists for every $T\in\intcurr_k(X)$ with $\bdry T=0$ an $S\in\intcurr_{k+1}(X)$ with $\bdry S=T$ and
 \begin{equation*}
  \mass{S}\leq C'[\mass{T}]^{\alpha}
 \end{equation*}
and which has the following property: Whenever $x\in\spt S$ and $0\leq s\leq \dist(x,\spt T)$ then
\begin{equation}\label{equation:basic-growth-estimate-filling}
 \|S\|(B(x,s))\geq \frac{s^{k+1}}{(3C)^k(k+1)^{k+1}}.
\end{equation}
Moreover, if $3C(k+1)\leq s\leq \dist(x,\spt T)$ then
\begin{equation*}
 \|S\|(B(x,s))\geq \left\{\begin{array}{ll}
  3C\left\{1+\frac{\alpha-1}{\alpha}\left[\frac{s}{3C}-(k+1)\right]\right\}^{\frac{\alpha}{\alpha-1}} &\text{if $\alpha>1$}\\*[0.5ex]
  3C\exp({\frac{s}{3C}-(k+1)}) & \text{if $\alpha =1$.}\\
\end{array}
\right.
\end{equation*}
\el
\br\label{remark:absolute-minimizers-growth-also}{\rm
It will be clear from the proof that all the conclusions hold also true for absolutely area minimizing currents $S\in\intcurr_{k+1}(X)$ with $\bdry S=T$.
}
\er
\begin{proof}
Let ${\mathcal M}$ denote the complete metric space consisting of all $S\in\intcurr_{k+1}(X)$ with $\bdry S=T$ and endowed with the 
metric given by $d_{\mathcal M}(S,S'):= \mass{S-S'}$. 
Choose an $\tilde{S}\in{\mathcal M}$ satisfying
$\mass{\tilde{S}}\leq C'[\mass{T}]^\alpha$. By the variation principle in \cite{Ekeland} there exists an $S\in{\mathcal M}$ with 
$\mass{S}\leq \mass{\tilde{S}}$ 
and such that the function 
\begin{equation*}
 S'\mapsto \mass{S'}+\frac{1}{2}\mass{S-S'}
\end{equation*}
has a minimum at $S'=S$. Let $x\in\spt S\ohne\spt T$ and set $\varrho_x(y):= d(x,y)$. 
Then the slicing theorem \cite[Theorems 5.6 and 5.7]{Ambr-Kirch-curr} implies that for almost every $0<s<\dist(x,\spt T)$
the slice $\slice{S}{\varrho_x}{s}$ exists, has zero boundary, and belongs to $\intcurr_{k}(X)$. For an
$S_s\in\intcurr_{k+1}(X)$ with $\bdry S_s=\slice{S}{\varrho_x}{s}$ the integral current
$S\rstr (X\ohne B(x,s))+S_s$ has boundary $T$ and thus, comparison with $S$ yields
\begin{equation*}
\mass{S\rstr (X\ohne B(x,s))+S_s} + \frac{1}{2}\mass{S\rstr B(x,s)-S_s}\geq \mass{S}.
\end{equation*}
If, moreover $S_s$ is chosen such that $\mass{S_s}\leq C[\mass{\slice{S}{\varrho_x}{s}}]^{\alpha}$ then the above estimate implies that
\begin{equation*}
\mass{S\rstr B(x,s)}\leq3\mass{S_s}\leq 3C[\mass{\slice{S}{\varrho_x}{s}}]^{\alpha}                                                   
\end{equation*}
for almost every $s\in(0,\dist(x,\spt T))$.
Since, by the slicing theorem, $\mass{\slice{S}{\varrho_x}{s}}\leq \beta'(s)$ for almost every $s\in(0,\dist(x,\spt T))$, where
$\beta(s):= \|S\|(B(x,s))$, we obtain
\begin{equation}\label{equation:diff-inequality-beta}
 \beta(s)\leq 3C[\beta'(s)]^{\alpha}
\end{equation}
for almost every $s\in(0,\dist(x,\spt T))$. 
Using the isoperimetric inequality of Euclidean type, i.e. with power $\frac{k+1}{k}$, proved in \cite[Theorem 1.2]{Wenger-GAFA}, the same arguments as above yield
\begin{equation*}
 \beta(s)\leq 3C\beta'(s)^{\frac{k+1}{k}}
\end{equation*}
from which follows
\begin{equation*}
 \|S\|(B(x,s))\geq \frac{s^{k+1}}{(3C)^k(k+1)^{k+1}}\quad\text{for all $0\leq s \leq\dist(x,\spt T).$}
\end{equation*}
This proves \eqref{equation:basic-growth-estimate-filling}. In particular, we have $\beta(3C(k+1))\geq 3C$. If $\alpha>1$ it follows from 
\eqref{equation:diff-inequality-beta} that
\begin{equation*}
 \begin{split}
 \beta(s)^{\frac{\alpha-1}{\alpha}}&\geq \beta(3C(k+1))^{\frac{\alpha-1}{\alpha}} + \frac{\alpha-1}{\alpha}\cdot\frac{s-3C(k+1)}{(3C)^{1/\alpha}}\\
  &\geq (3C)^{\frac{\alpha-1}{\alpha}}\left\{1+\frac{\alpha-1}{\alpha}\left[\frac{s}{3C}-(k+1)\right]\right\}
 \end{split}
\end{equation*}
from which the second statement follows for $\alpha>1$. On the other hand, if $\alpha=1$ then
\begin{equation*}
 \beta(s)\leq 3C\beta'(s)
\end{equation*}
and hence 
\begin{equation*}
 \frac{s-3C(k+1)}{3C}\leq \ln\left(\frac{\beta(s)}{\beta(3C(k+1))}\right)\leq \ln\left(\frac{\beta(s)}{3C}\right).
\end{equation*}
This proves the lemma.
\end{proof}
We are now ready to prove the main theorem. The proof of the first part of the theorem is a variation of the arguments given in \cite{Lang-Schroeder}. 
\begin{proof}[{Proof of \thmref{Theorem:rank-div}}]
 We begin by proving the first statement of the theorem. For this let $F\subset X$ be a flat of maximal dimension $k+1=\Erank X$. 
 Choose $x_0\in F$, set $A_0:=\hm^k(S(x_0, 1)\cap F)$, and let $\varrho\in(0,1)$ be arbitrary.
 We will show that 
 \begin{equation}\label{equation:rank-div-first-part-claim}
  \liminf_{r\to\infty}\frac{\delta^k_{x_0, \varrho, A_0}}{r^{k+2}}>0.
 \end{equation}
 By \lemref{Lemma:basepoint-delta} this will imply $\divGMT_k(X)\succeq r^{k+2}$. In order to prove \eqref{equation:rank-div-first-part-claim} fix
 $r> 24\sqrt{k+1}s_0/\varrho$ where $s_0>0$ is chosen large enough as below and set $B:= B(x_0, r)\cap F$. Denote by $\iota: B\hookrightarrow X$ the inclusion
 map and set $T:= \bdry(\iota_{\#}\Lbrack \chi_B\Rbrack)$. Note that $T$ is an element of $\intcurr_k(X)$ and  
 satisfies $\mass{T}=\hm^k(S(x_0, r)\cap F)=A_0r^k$. Let furthermore $S\in\intcurr_{k+1}(X)$ be such that $\bdry S= T$ and $\spt S \subset X\ohne U(x_0, \varrho r)$.
Denote by $\pi: X\to F$ the orthogonal projection onto $F$ and observe that $\pi_{\#}S=\iota_{\#}\Lbrack \chi_B\Rbrack$
by the constancy theorem \cite[4.1.7]{Federer}.
Let now $Q\subset F\cap B(x_0,\varrho r/4)$ be a closed $(k+1)$-dimensional Euclidean cube of edge length $3s_0$. 
We adapt the argument in the proof of Proposition 3.2 in \cite{Lang-Schroeder} to estimate $\|S\|(\pi^{-1}(Q))$ from below. For this set $\nu:=2^{k+1}$ and
let the vertices $q_1,\dots,q_{\nu}$ of $Q$ be ordered in such a way that each segment $[q_i,q_{i+1}]$, $i=1,\dots, \nu-1$, is an edge of $Q$. 
Set $P=\cup_{i=1}^{\nu - 1}[q_i,q_{i+1}]$ and denote by $R$ the union of all $(k+1)$-cubes of $Q$ with edge length $s_0$ which
meet $P$. Let $Q_i\subset Q$ be the cube of edge length $s_0$ which contains $q_i$. Denoting by $Z$ the common $k$-face of $R$ and $Q_1$ 
we define a $1$-Lipschitz map $\psi: R\to Z$ in such a way that each fiber $\psi^{-1}(\{z\})$ is a connected polygonal arc lying at constant distance from $P$.
After possibly changing the edge length of $Q$ by an arbitrarily small amount we may assume by the slicing theorem that 
$S\rstr\pi^{-1}(Q)\in\intcurr_{k+1}(X)$. Furthermore, 
$S_z:=\langle S\rstr\pi^{-1}(Q), \psi\circ\pi, z\rangle\in\intcurr_1(X)$ for almost every $z\in Z$. 
Since 
\begin{equation*}
\spt \bdry(S\rstr \pi^{-1}(Q))\subset \pi^{-1}\left(\overline{F\ohne Q}\right)
\end{equation*}
and
$\bdry S_z = (-1)^k\langle \bdry(S\rstr\pi^{-1}(Q)), \psi\circ\pi, z\rangle$ it follows that $\bdry S_z$ is supported in $\pi^{-1}\left(\overline{F\ohne Q}\right)$.
Furthermore we have 
\begin{equation*}
 \pi_{\#}S_z = \langle\pi_{\#}(S\rstr\pi^{-1}(Q)),\psi, z\rangle=\langle(\pi_{\#}S)\rstr Q, \psi,z\rangle=\Lbrack\psi^{-1}(\{z\})\Rbrack
\end{equation*}
which implies that $\spt S_z$ has a connected component whose image under $\pi$ is $\psi^{-1}(\{z\})$.
Therefore, if $s_0$ is chosen suitably large (only depending on $X$) as in Lemma 3.1 in \cite{Lang-Schroeder} it follows that
\begin{equation*}
 \mass{S_z}\geq  \lambda \dist(\spt S_z, F)\geq\lambda \varrho r/4
\end{equation*}
where $\lambda>0$ is a constant only depending on $X$. Application of the slicing theorem finally yields
\begin{equation*}
 \|S\|(\pi^{-1}(Q))\geq \int_Z \mass{S_z}dz \geq \frac{\lambda\varrho}{4}\; s_0^k r.
\end{equation*}
Since $B(x_0,r\varrho/4)\cap F$ contains at least 
\begin{equation*}
\left (\frac{\varrho r}{6\sqrt{k+1}s_0}-1\right)^{k+1}
\end{equation*}
pairwise disjoint cubes of edge length $3s_0$ we obtain
\begin{equation*}
 \mass{S}\geq Er^{k+2}
\end{equation*}
for a constant $E$ depending only on $\varrho$, $\lambda$, $s_0$, and $k$. This completes the proof of the first statement.

We now turn to the proof of the second part of the theorem. The proof is by contradiction and we therefore 
assume that $\divGMT_k(X)$ grows faster than $r^{k+1}$. There thus exist for $\varrho:=\frac{1}{4}$ an $x_0\in X$ and $A>0$ such that
 \begin{equation*}
  \limsup_{r\to\infty}\frac{\delta^k_{x_0,\varrho,A}(r)}{r^{k+1}}=\infty.
 \end{equation*}
 In particular, there exist a sequence $r_m\nearrow\infty$ and a 
 sequence $(T_m)\subset\intcurr_k(X)$ with $\bdry T_m=0$ and such that furthermore:
 \begin{enumerate}
  \item $\spt T_m\subset S(x_0, r_m)$ compact for all $m\in\N$
  \item $r_m^{-(k+1)}\fillvol_{X\ohne U(x_0,r_m/4)}(T_m)\to\infty$ as $m\to\infty$
  \item $\mass{T_m}\leq A r_m^k$ for every $m\in\N$.
 \end{enumerate}
By the compactness and closure theorems \cite[Theorems 5.2 and 8.5]{Ambr-Kirch-curr} 
there exists for every $m\in\N$ an absolutely area minimizing $S_m\in\intcurr_{k+1}(X)$ 
with $\bdry S_m=T_m$ and with compact support. By \cite[Theorem 1.2]{Wenger-GAFA} we have
 \begin{equation}\label{equation:mass-bound-filling}
  \mass{S_m}\leq C[\mass{T_m}]^{\frac{k+1}{k}}\leq CA^{\frac{k+1}{k}}r_m^{k+1}
 \end{equation}
for some constant $C$ depending only on $k$, and in view of property (ii) we may therefore assume that
\begin{equation*}
 \spt S_m\cap B(x_0, r_m/4)\not=\emptyset
\end{equation*}
for every $m\in\N$.
Fix a point $x_m$ in the intersection. By the slicing theorem and inequality
\eqref{equation:mass-bound-filling} there exists $r'_m\in(\frac{1}{2}r_m,\frac{3}{4}r_m)$ such that $S_m\rstr B(x_0,r'_m)\in\intcurr_{k+1}(X)$
and
\begin{equation*}
 \mass{\bdry(S_m\rstr B(x_0,r'_m))}\leq 4CA^{\frac{k+1}{k}}r_m^k.
\end{equation*}
Define a metric space $(Y,d_Y)$ as the disjoint union $Y:=\bigsqcup_{m=1}^\infty Y_m$ where $Y_m:= B(x_0, r_m)$ and with the metric on $Y$ given by
\begin{equation*}
 d_Y(y,y'):=\left\{
 \begin{array}{c@{\quad}l}
  \frac{1}{r_m}d(y,y') &\text{if $y,y'\in Y_m$ for some $m\in\N$}\\
  3 &\text{otherwise.}
 \end{array}\right.
\end{equation*}
Note that balls of radius strictly less than $3$ in $Y$ are ${\rm CAT}(0)$ and hence are geodesic and admit cone type inequalities as in 
\eqref{equation:non-optimal-cone-type-cat0}. These two properties will be needed to invoke \cite[Theorem 1.4]{Wenger-flatconv} at a later stage.
Denote by $S'_m$ the current $S_m\rstr B(x_0, r'_m)$ viewed as an element of $\intcurr_{k+1}(Y)$ with 
\begin{equation*}
 Z_m:= \spt S'_m\subset Y_m\subset Y
\end{equation*}
and observe that $\mass{S'_m}\leq CA^{\frac{k+1}{k}}$ and $\mass{\bdry S'_m}\leq 4CA^{\frac{k+1}{k}}$.
It follows from \eqref{equation:basic-growth-estimate-filling} that the sequence $(Z_m,d_Y)$ is equi-compact and equi-bounded and hence,
by Gromov's compactness theorem, there exists (after passage to a subsequence) a compact metric space $(Z,d_Z)$ and isometric embeddings 
$\varphi_m: (Z_m, d_Y)\hookrightarrow (Z,d_Z)$. Furthermore, we may assume without loss of generality 
that $\varphi_m(Z_m)$ converges to to some compact subset $Z'\subset Z$ with respect to the Hausdorff metric and, by the compactness and closure theorems for currents,
that $\varphi_{m\#}S'_m$ weakly converges to some $S\in\intcurr_{k+1}(Z)$. 

We now claim that $S\not=0$. Before proving this claim we show how the theorem follows from it.
For this, fix an ultrafilter $\omega$ on $\N$ and denote by $X_\omega$ the asymptotic cone associated with the sequence 
$(X,\frac{1}{r_m}d_X, x_0)$. We construct a map $\psi:Z'\to X_\omega$ as follows.
For $z\in Z'$ there exists $z_m\in Z_m$ such that $\varphi_m(z_m)\to z$. We set $\psi(z):= {(z_m)}_{m\in\N}$. It is straight forward to check that
$\psi$ is well-defined and an isometric embedding. Since $\spt S\subset Z'$ and since $S\not=0$ we obtain that $\psi_{\#}S$ is a non-zero $(k+1)$-dimensional 
integral current in $X_\omega$. By the representation formula for integer rectifiable currents (see \secref{Section:Preliminaries} or Theorem 4.5 in 
\cite{Ambr-Kirch-curr}) there then exists a bi-Lipschitz map $\nu:K\subset\R^{k+1}\to X_\omega$ where
$K$ is measurable and of strictly positive Lebesgue measure. From \eqref{equation:metric-differential-strong} and \cite[Theorem A]{Kleiner-local-structure} 
we conclude that the geometric dimension of $X_\omega$ is at least $k+1$ and hence, by Theorem C of \cite{Kleiner-local-structure}, that the Euclidean rank of $X$ is at 
least $k+1$. This contradicts our assumption that $\Erank X\leq k$ and concludes the proof of the theorem under the assumption that the claim holds.

We return to the proof of the claim and assume first that $\bdry S'_m$ converges weakly to $0$. By Theorem 1.4 of \cite{Wenger-flatconv} we then have
$\fillvol(\bdry S'_m)\to 0$. In particular, there exist absolutely area minimizing currents $S''_m\in\intcurr_{k+1}(Y)$ with $\bdry S''_m=\bdry S'_m$ and 
$\spt S''_m\subset Y_m$ for all $m\in\N$ and
such that $\mass{S''_m}\to 0$. Denote by $\tilde{S}_m$ the current $S''_m$ viewed as an integral current in $X$. Then $\tilde{S}_m$ is absolutely area minimizing
and satisfies $\bdry\tilde{S}_m=\bdry(S_m\rstr B(x_0, r'_m))$ and
\begin{equation*}
 \frac{\mass{\tilde{S}_m}}{r_m^{k+1}}=\mass{S''_m}\to 0.
\end{equation*}
It follows from \eqref{equation:basic-growth-estimate-filling} that 
\begin{equation*}
\spt\tilde{S}_m\subset X\ohne U(x_0,r_m/4)
\end{equation*}
for $m$ large enough. This leads to a contradiction with (ii). Indeed, $S_m\rstr (X\ohne B(x_0, r'_m))+\tilde{S}_m$ is a filling of $T_m$ with support outside 
$U(x_0, r_m/4)$ and with mass bounded from above by $Dr_m^{k+1}$ for a suitable constant $D$.
This shows that $\bdry S'_m$ does not weakly converge to $0$. In particular, there exist Lipschitz maps 
$f,\pi_1,\dots,\pi_k\in\lip(Y)$ and $\varepsilon>0$ such that (after passage to a subsequence)
\begin{equation*}
 \bdry S'_m(f,\pi_1,\dots,\pi_k)\geq \varepsilon\quad\text{for all $m\in\N$.}
\end{equation*}
Note that since $Y$ is a bounded metric space the functions $f$ and $\pi_i$ are bounded. We define Lipschitz functions $f_m$ and $\pi_i^m$ on $\varphi_m(Z_m)$ by
$f_m(z):= f(\varphi_m^{-1}(z))$ and $\pi_i^m(z):=\pi_i(\varphi_m^{-1}(z))$ for $z\in\varphi_m(Z_m)$. Here, we view $\varphi^{-1}_m$ as a map from $\varphi(Z_m)$ 
to $Y$ with image in $Y_m\subset Y$. 
By McShane's Lipschitz extension theorem there exist extensions $\hat{f}_m, \hat{\pi}_i^m: Z\to\R$ of $f_m$ and $\pi_i^m$ with the same Lipschitz constants
as $f$ and $\pi_i$, which we may assume to be $1$.
By Arzel\`a-Ascoli theorem we may assume that $\hat{f}_m$ and 
$\hat{\pi}_i^m$ converge uniformly to Lipschitz maps $\hat{f}$, $\hat{\pi}_i$ on $Z$. 
Finally we abbreviate $T'_m:=\varphi_{m\#}(\bdry S'_m)$ and use 
\cite[Proposition 5.1]{Ambr-Kirch-curr} to estimate
\begin{align*}
   \bdry S (\hat{f},\hat{\pi}_1,\dots,\hat{\pi}_k) &= \lim_{m\to\infty}T'_m(\hat{f},\hat{\pi}_1,\dots,\hat{\pi}_k)\\
&=\lim_{m\to\infty}\Big[T'_m(\hat{f}_m,\hat{\pi}_1^m,\dots,\hat{\pi}_k^m) + T'_m(\hat{f}-\hat{f}_m,\hat{\pi}_1,\dots,\hat{\pi}_k)\\
  &\quad\qquad\qquad+T'_m(\hat{f}_m,\hat{\pi}_1,\dots,\hat{\pi}_k)-T'_m(\hat{f}_m,\hat{\pi}_1^m,\dots,\hat{\pi}_k^m)\Big]\\
  &\geq \varepsilon - \limsup_{m\to\infty}\Big[\prod_{i=1}^k\lip(\hat{\pi}_i)\int_Z|\hat{f}-\hat{f}_m|\,d\|T'_m\|\Big]\\
  & \qquad- \limsup_{m\to\infty}\Big[\lip(\hat{f}_m)\sum_{i=1}^k\int_Z|\hat{\pi}_i-\hat{\pi}_i^m|\,d\|T'_m\|\Big]\\
  &= \varepsilon.
\end{align*}
This shows that indeed $\bdry S\not=0$ and hence also $S\not=0$ and therefore completes the proof of the claim and of the theorem.

\end{proof}
\br
{\rm
 We point out that for the second part of the statement properness of $X$ is only needed for Theorem C of \cite{Kleiner-local-structure}. The existence of the absolutely
 area minimizing currents $S_m$ when $X$ is not proper follows from \cite[Theorem 1.6]{Wenger-GAFA}.
}
\er

We mention that if $X$ is a cocompact Hadamard manifold one can give a proof of the second part of the theorem without using asymptotic cones. We can use 
(the proof of) Theorem~1 in \cite{Anderson-Schroeder} in the following way instead. Let $S_m$ and $x_m$ be as in our proof. By cocompactness and the
compactness theorem for integral currents we may assume without loss of generality that $x_m$ converges to some $y\in X$ and that $S_m$ converges to 
an absolutely area minimizing local integral current $\Sigma\in\intcurr_{k+1,\;\operatorname{loc}}(X)$ with $\bdry\Sigma=0$. Using the monotonicity formula, inequality \eqref{equation:mass-bound-filling} and the fact that $\dist(x_m, \spt \bdry S_m)\geq \frac{3r_m}{4}$ one readily obtains
\begin{equation*}
 \omega_{k+1}r^{k+1}\leq \|\Sigma\|(B(y, r))\leq Cr^{k+1}
\end{equation*}
for some constant $C$ and for all $r>0$. Then the proof of Theorem~1 in \cite{Anderson-Schroeder} (see (3.3) and thereafter) yields the existence of a 
$(k+1)$-dimensional flat in $X$ which contradicts the assumption that $\Erank X\leq k$.

\section{Divergence versus isoperimetric inequality}\label{Section:div-isop}
We first prove \propref{Proposition:isop-div}.
\begin{proof}[{Proof of \propref{Proposition:isop-div}}]
 Fix $x_0\in X$ and $A>0$, set $\varrho_0:=\frac{1}{2}$ and let
 $r>0$ be large enough (as chosen below). For a $T\in\intcurr_k(X)$ with $\spt T\subset S(x_0, r)$, $\bdry T=0$ and $\mass{T}\leq Ar^k$ let $S\in\intcurr_{k+1}(X)$
 be as in \lemref{Lemma:Q-min-growth}. We show that $S$ is an admissible filling in the sense that its support is contained in $X\ohne U(x_0,r/2)$. 
 Indeed, we have 
 \begin{equation*}
  \mass{S}\leq C'[\mass{T}]^\alpha\leq C'A^\alpha r^{k\alpha}
 \end{equation*}
which together with the growth estimate in \lemref{Lemma:Q-min-growth} yields for $x\in\spt S$
\begin{equation*}
 \dist(x,\spt T)\leq \left\{\begin{array}{ll}
   D_1+D_2r^{k(\alpha-1)} & \text{if $\alpha>1$}\\
   D_1+D_2\ln r & \text{if $\alpha=1$}
  \end{array}\right.
\end{equation*}
where $D_1,D_2$ are constants which only depend on $C$, $k$, $\alpha$, and $A$.
Here, $C$ and $C'$ are the constants from \lemref{Lemma:Q-min-growth}.
Since $k(\alpha - 1)<1$ it follows that for $r>0$ large enough $\spt S$ lies outside the ball $B(x_0, r/2)$. This concludes the proof.
\end{proof}

\subsection{Optimal cone inequality and linear isoperimetric inequality }
In this section we prove an optimal cone inequality for integral currents in ${\rm CAT}(\kappa)$-spaces. The proof is a generalization of 
the proof in the setting of manifolds. However, since we can only work with comparison triangles and not with Jacobi fields the estimates become more elaborate and
we need to use the metric derivative of Lipschitz maps into metric spaces as well.
As a direct consequence of the cone inequality we will obtain \thmref{Theorem:linear-isoperimetric-inequality} as well as the monotonicity formula for absolutely 
area minimizing currents in Hadamard spaces. All of these results are well-known in the smooth setting, i.e.~in Hadamard manifolds.

It is not difficult to prove (see \cite[Proposition 2.10]{Wenger-GAFA}) that Hadamard spaces $X$ admit cone type inequalities for $\intcurr_k(X)$, $k\geq 0$, 
in the following sense: If $T\in\intcurr_k(X)$ satisfies $\bdry T=0$, or $T(1)=0$ if $k=0$, and has bounded support then there exists an $S\in\intcurr_{k+1}(X)$ with 
$\bdry S=T$ and
\begin{equation}\label{equation:non-optimal-cone-type-cat0}
 \mass{S}\leq (k+1)\diam(\spt T)\mass{T}.
\end{equation}
The constant $(k+1)$ appearing in \eqref{equation:non-optimal-cone-type-cat0} is not optimal.
In order to state the optimal cone inequality let
$(X,d)$ be a complete ${\rm CAT}(\kappa)$-space, $\kappa\in\R$, and fix a point $x_0\in X$. For $x\in X$ with $d(x,x_0)<\frac{1}{2}D_\kappa$ 
we set $\varphi(t,x):=c_{x_0x}(t)$ where $c_{x_0x}:[0,1]\to X$ denotes the constant speed geodesic from $x_0$ to $x$. We furthermore define
\begin{equation*}
 s_\kappa(r):=\left\{\begin{array}{c@{\qquad}l} \frac{1}{\sqrt{-\kappa}}\sinh(\sqrt{-\kappa}r) &\text{if $\kappa<0$}\\
                                               r        &\text{if $\kappa=0$}\\
                                               \frac{1}{\sqrt{\kappa}}\sin(\sqrt{\kappa}r)  &\text{if $\kappa>0$}
                     \end{array}\right.
\end{equation*}
for $r\geq 0$.
We note that $s_\kappa$ is the norm of a normal Jacobi field $Y$ along a geodesic parameterized by arc-length
on $\spaceform_\kappa^2$ with $Y(0)=0$ and $\|Y'(0)\|=1$.
\bt\label{theorem:cat-cone-inequality}
 Let $(X,d)$ be a complete ${\rm CAT}(\kappa)$-space and $x_0\in X$. Let $k\geq 0$ and suppose
 $T\in\intcurr_k(X)$ has support in the ball $B(x_0,R)$, where $R<\frac{1}{2}D_\kappa$ in case $\kappa>0$.
 Then $\cone{x_0}{T}:= \varphi_\#([0,1]\times T)$ satisfies
\begin{equation*}
 \mass{\cone{x_0}{T}}\leq \int_{B(x_0,R)}\frac{d(x_0,x)}{[s_\kappa(d(x_0,x))]^k}\int_0^1[s_\kappa(td(x_0,x))]^kdt\; d\|T\|(x).
\end{equation*}
In particular, if $\kappa=0$ then we obtain
\begin{equation}\label{equation:optimal-cone}
 \mass{\cone{x_0}{T}}\leq \frac{R}{k+1}\mass{T}.
\end{equation}
\et
\br{\rm
For Banach spaces inequality \eqref{equation:non-optimal-cone-type-cat0} still holds, whereas \eqref{equation:optimal-cone} is in general false 
(even for $2$-dimensional normed spaces) as Example 10.3 in \cite{Ambr-Kirch-curr} illustrates.
}
\er
\bl\label{lemma:alexandrov-metric-differential}
 Let $(X,d)$ be an Alexandrov space of curvature bounded above and $\varphi:K\to X$ a Lipschitz map where $K\subset\R^k$ is measurable and $k\geq 1$. 
 Suppose $z\in K$ is a density point, $\md\varphi_z$ exists, is non-degenerate and satisfies \eqref{equation:metric-differential-strong}. Then 
 $\md\varphi_z$ is induced by an inner product.
\el
\begin{proof}
 If $k=1$ then the statement holds for trivial reasons. We may therefore assume that $k\geq 2$.
 Let $z\in K$ be as in the assumption. We show that $(\R^k, \md\varphi_z)$ is an inner product space. By \cite[Proposition II.1.11]{Bridson-Haefliger}
 it is enough to prove
 that for any four vectors $v_1,\dots,v_4\in\R^k$ there exists a subembedding in $\R^2$, i.e.~there exist four points 
 $w_1,\dots, w_4\in\R^2$ such that $\md\varphi_z(v_i-v_{i+1})= |w_i-w_{i+1}|$, $i\in\N$ mod $4$, as well as
 $\md\varphi_z(v_1-v_3)\leq |w_1-w_3|$ and $\md\varphi_z(v_2-v_4)\leq |w_2-w_4|$.
 Since every normed space which is ${\rm CAT}(0)$ is in fact a pre-Hilbert space this will prove the lemma.\\
 To prove the existence of a subembedding let
 $\varepsilon>0$ be such that $B:=B(\varphi(x),\varepsilon)$ is ${\rm CAT}(\kappa)$ for some $\kappa\in\R$. Then $B$ 
 endowed with the rescaled metric $d_r:=\frac{1}{r}d$ is ${\rm CAT}(\kappa\sqrt{r})$ for every $r>0$.
 We first assume $\kappa\leq 0$ so that $(B,d_r)$ is ${\rm CAT}(0)$.
 Since $z$ is a Lebesgue density point of $K$ we may assume after approximation that $z+rv_i\in K$ for $r>0$ small enough and for $i=1,\dots,4$. For $r>0$ 
 sufficiently small there exist again by \cite[Proposition II.1.11]{Bridson-Haefliger} points $w_1^r,\dots,w_4^r\in\R^2$ such that
 \begin{itemize}
  \item $\frac{1}{r}d(\varphi(z+rv_i),\varphi(z+rv_{i+1}))= |w_i^r-w_{i+1}^r|$ for $i\in\N$ mod $4$
  \item $\frac{1}{r}d(\varphi(z+rv_1),\varphi(z+rv_3))\leq |w_1^r-w_3^r|$
  \item $\frac{1}{r}d(\varphi(z+rv_2),\varphi(z+rv_4))\leq |w_2^r-w_4^r|$.
 \end{itemize}
 Of course it is not restrictive to assume that $w_1^r=0$ for all $r>0$. By the Lipschitz continuity of $\varphi$ we conclude that all $w_i^r$ lie in a fixed ball
 centered at $0$.
 There then exists a sequence $r_n$ converging to $0$ such that $w_i^{r_n}$ converges to some $w_i$ as $n\to\infty$ for every $i$. It then follows immediately 
from property \eqref{equation:metric-differential-strong} that $w_1,\dots,w_4$ constitute a subembedding of $v_1,\dots,v_4$. 
 This concludes the proof in the case $\kappa\leq 0$. The
 case $k>0$ is almost analogous. It is enough to note that the comparison spaces $\spaceform^2_{\kappa\sqrt{r}}$ in which the comparison points $w^r_i$ lie
 converge to $\R^2$ in the pointed Hausdorff-Gromov metric (the base point being the north pole) as $r\searrow 0$.
\end{proof}
\begin{proof}[{Proof of \thmref{theorem:cat-cone-inequality}}]
 If $k=0$ then we may assume without loss of generality that $T$ is of the form
 \begin{equation*}
  T(f)=\theta_1 f(x_1),\quad\text{$f$ bounded and Lipschitz,}
 \end{equation*}
for some $x_1\in X$ and $\theta_1\in\Z\ohne\{0\}$. Then
\begin{equation*}
 \begin{split}
  (\cone{x_0}{T})(f,\pi)&=([0,1]\times T)(f\circ\varphi,\pi\circ\varphi)\\
   &=\int_0^1T\left(f\circ\varphi(t,\cdot)\frac{\partial (\pi\circ\varphi)}{\partial t}\right)dt\\
   &=\theta_1\int_0^1f(\varphi(t,x_1))\frac{\partial(\pi\circ\varphi(t,x_1))}{\partial t}dt
 \end{split}
\end{equation*}
 from which we easily infer that
\begin{equation*}
 \mass{\cone{x_0}{T}}\leq |\theta_1|d(x_0,x_1)\leq R\mass{T}.
\end{equation*}
 This concludes the proof of the case $k=0$.

 If $k\geq 1$ then, by Theorem 4.5 in \cite{Ambr-Kirch-curr}, it is not restrictive to assume that $T=\psi_\#\Lbrack\theta\Rbrack$ for some bi-Lipschitz map 
 $\psi:K\to X$, $K\subset\R^k$ compact, and $\theta\in L^1(K,\Z)$. We give an explicit formula for $\cone{x_0}{T}$. For this, let $(f,\pi_1,\dots,\pi_k)\in\form^k(X)$
 with $\lip(\pi_i)\leq 1$, $i=1,\dots, k$. 
 We define $g_t:=g\circ\varphi(t,\cdot)$ whenever $g\in\lip(X)$ and furthermore write $\hat{\pi}^i:= (\pi_1,\dots,\pi_{i-1},\pi_{i+1},\dots,\pi_{k+1})$. 
 We compute
 \begin{equation*}
  \begin{split}
   (\cone{x_0}{T}&)(f,\pi_1,\dots,\pi_{k+1})\\
    &= ([0,1]\times T)(f\circ\varphi,\pi_1\circ\varphi,\dots,\pi_{k+1}\circ\varphi)\\
    &=\sum_{i=1}^{k+1}(-1)^{i+1}\int_0^1 T\left(f_t\frac{\partial \pi_{i\,t}}{\partial t},\pi_{1\,t},
                                               \dots,\pi_{i-1\,t},\pi_{i+1\,t},\dots,\pi_{k+1\,t}\right)dt\\
    &=\sum_{i=1}^{k+1}(-1)^{i+1}\int_0^1\int_K\theta\; f\circ\tilde{\varphi}\;\frac{\partial(\pi_i\circ\tilde{\varphi})}{\partial t}
                                                \det(\nabla_{\R^k}(\hat{\pi}^i\circ\tilde{\varphi}))\,d\lm^k\,dt\\
    &= \int_{[0,1]\times K} \tilde{\theta}\; f\circ\tilde{\varphi}\; \det(\nabla(\pi\circ\tilde{\varphi}))\,d\lm^{k+1}\\
    &=(\tilde{\varphi}_\#\Lbrack\tilde{\theta}\Rbrack)(f,\pi_1,\dots,\pi_{k+1}),
  \end{split}
 \end{equation*}
where we have furthermore set $\tilde{\varphi}(t,z):= \varphi(t,\psi(z))$ and $\tilde{\theta}(t,z):= \theta(z).$
Since 
\begin{equation*}
 \left|\det\left(\frac{\partial(\pi_i\circ\tilde{\varphi})(t,z)}{\partial x_j}\right)\right|\leq \lambda_{\md\tilde{\varphi}_{(t,z)}}
  \jac_{k+1}(\md\tilde{\varphi}_{(t,z)})
\end{equation*}
we easily obtain using the definition of mass and the area formula for Lipschitz maps \cite[Corollary 8]{Kirchheim} that
\begin{equation}\label{equation:mass-representation-param}
  \mass{\tilde{\varphi}_\#\Lbrack\tilde{\theta}\Rbrack}\leq \int_A |\theta(z)|\lambda_{\md\tilde{\varphi}_{(t,z)}}
                                                                                \jac_{k+1}(\md\tilde{\varphi}_{(t,z)})\,d\lm^{k+1}(t,z).
 \end{equation}
where $A$ is the set of points $(t,z)\in[0,1]\times K$ such that $z$ is a Lebesgue point of $K$ and 
$\md\tilde{\varphi}_{(t,z)}$ exists, is non-degenerate, and satisfies 
\eqref{equation:metric-differential-strong}.
By \lemref{lemma:alexandrov-metric-differential} above $\md\tilde{\varphi}_{(t,z)}$ comes from an inner product for all $(t,z)\in A$ and therefore
$\lambda_{\md\tilde{\varphi}(t,z)}=1$. 
We now estimate $\jac_{k+1}(\md\tilde{\varphi}_{(t,z)})$. 
We may assume without loss of generality that $\md\psi_z$ exists and is induced by an inner product and that the 
function $\tilde{\nu}(z'):= d(x_0,\psi(z'))$ is differentiable at the point $z$. Choose a basis 
$\{v_1,\dots,v_k\}$ of $\R^k$ which is orthonormal with respect to the inner product inducing $\md\psi_z$ and for which $v_1$ is parallel to 
$\nabla \tilde{\nu}(z)$. If $\nabla\tilde{\nu}(z)=0$ then we do not pose any restriction on the choice of $v_1$.
Set 
\begin{equation*}
 Q:= \left\{\left(s,\sum_{i=1}^k r_iv_i\right): 0\leq s\leq 1,0\leq r_i\leq 1\right\}\subset \R\times\R^k
\end{equation*}
and observe that
\begin{equation*}
  \jac_{k+1}(\md\tilde{\varphi}_{(t,z)}) = 
     \frac{\hm^{k+1}_{\md\tilde{\varphi}_{(t,z)}}(Q)}{\lm^{k+1}(Q)}= \hm^{k+1}_{\md\tilde{\varphi}_{(t,z)}}(Q)\jac_k(\md\psi_z)
\end{equation*}
and
\begin{equation*}
 \hm^{k+1}_{\md\tilde{\varphi}_{(t,z)}}(Q)\leq \md\tilde{\varphi}_{(t,z)}(1,0)\md\tilde{\varphi}_{(t,z)}(h,v_1)\prod_{i=2}^k\md\tilde{\varphi}_{(t,z)}(0,v_i)
\end{equation*}
for every $h\in\R$. The latter is a consequence of the fact that $\md\tilde{\varphi}_{(t,z)}$ comes from an inner product.

We first estimate $\md\tilde{\varphi}_{(t,z)}(0,v_i)$ from above for $i\geq 2$. Fix an $i\geq 2$ and set $v:=v_i$. For $r>0$ small enough let
$\overline{\Delta}(\overline{x}_0,\overline{x},\overline{x}_r)$ be the comparison triangle in $\spaceform_\kappa^2$ of the triangle 
$\Delta(x_0,\psi(z),\psi(z+rv))$. By a simple approximation argument we may assume that $z+rv\in K$ for all $r>0$ small enough.
Denoting by $\alpha_r$ the angle at $\overline{x}$ between the geodesics $[\overline{x},\overline{x}_0]$ and 
$[\overline{x},\overline{x}_r]$ one verifies, using the law of cosines for $\spaceform^2_\kappa$, see \cite[I.2.13]{Bridson-Haefliger}, 
that $\alpha:= \lim_{r\searrow 0}\alpha_r$ exists and satisfies
\begin{equation*}
 \cos\alpha = -\frac{d\tilde{\nu}_z(v)}{\md\psi_z(v)}=0
\end{equation*}
and hence $\alpha=\frac{\pi}{2}$. 
Denote by $\overline{c}(\cdot,r)$ the geodesic parameterized on $[0,1]$ from $\overline{x}_0$ to $\overline{x}_r$. Then, by the ${\rm CAT}(\kappa)$-condition we obtain
\begin{equation*}
 \md\tilde{\varphi}_{(t,z)}(0,v) \leq \limsup_{r\searrow 0}\frac{1}{r}d_{\spaceform_\kappa^2}(\overline{c}(t,0),\overline{c}(t,r))= \|Y(t)\|,
\end{equation*}
where $Y$ is the normal Jacobi field along $\overline{c}(\cdot,0)$ with $Y(0)=0$ and $\|Y(1)\|=\md\psi_z(v)=1$, and consequently
\begin{equation*}
 \md\tilde{\varphi}_{(t,z)}(0,v_i)\leq\frac{s_\kappa(td(x_0,\psi(z)))}{s_\kappa(d(x_0,\psi(z)))}
\end{equation*}
for every $i\in\{2,3,\dots,k\}$.\\
Next, we estimate $\md\tilde{\varphi}_{(t,z)}(h,v_1)$ for a suitable $h\in\R$ by proceeding in a similar way.
Let $\overline{\Delta}(\overline{x}_0,\overline{x},\overline{x}_r)$, $\alpha_r$, and $\overline{c}$ be as above, but with $v:= v_1$. 
After possibly replacing $v_1$ by $-v_1$ we may assume that $\alpha_r\geq\frac{\pi}{2}$. Let $h(r)\in\R$ be such that the triangle 
$\Delta(\overline{x}_0,\overline{x},\overline{c}(1+rt^{-1}h(r),r))$ in $\spaceform_\kappa^2$ has a right angle at the vertex $\overline{x}$. Then one easily 
checks that
\begin{equation*}
 h(r)\to h:= -\frac{td\tilde{\nu}_z(v_1)}{\tilde{\nu}(z)}\leq0\qquad\text{as $r\searrow 0$}
\end{equation*}
and, since $d(\overline{x}_0,\overline{x}_r)\leq R<\frac{1}{2}D_\kappa$, 
\begin{equation*}
 \limsup_{r\searrow 0}\frac{1}{r}d_{\spaceform_\kappa^2}(\overline{x},\overline{c}(1+rt^{-1}h(r),r))
                               \leq\limsup_{r\searrow 0}\frac{1}{r}d_{\spaceform_\kappa^2}(\overline{x},\overline{x}_r)
\end{equation*}
and hence 
\begin{equation}\label{equation:jacobi-field-estimate}
 \limsup_{r\searrow 0}\frac{1}{r}d_{\spaceform_\kappa^2}(\overline{x},\overline{c}(1+rt^{-1}h(r),r))\leq\md\psi_z(v_1)=1.
\end{equation}
The ${\rm CAT}(\kappa)$-condition then implies
\begin{equation*}
 \begin{split}
 \md\tilde{\varphi}_{(t,z)}(h,v_1)&=\lim_{r\searrow 0}\frac{1}{r}d(\tilde{\varphi}(t+rh,z+rv_1),\tilde{\varphi}(t,z))\\
                                  &=\lim_{r\searrow 0}\frac{1}{r}d(\tilde{\varphi}(t(1+rt^{-1}h(r)),z+rv_1),\tilde{\varphi}(t,z))\\
                                  &\leq\limsup_{r\searrow 0}\frac{1}{r}d_{\spaceform_\kappa^2}(\overline{c}(t(1+rt^{-1}h(r)),r),\overline{c}(t,0))\\
                                  &=\|Y(t)\|,
 \end{split}
\end{equation*}
where $Y$ denotes the normal Jacobi field along $\overline{c}(\cdot,0)$ with $Y(0)=0$ and 
$$\|Y(1)\|=\limsup_{r\searrow 0}\frac{1}{r}d_{\spaceform_\kappa^2}(\overline{c}(1+rt^{-1}h(r),r),\overline{x}).$$ 
Together with \eqref{equation:jacobi-field-estimate} we obtain
\begin{equation*}
 \md\tilde{\varphi}_{(t,z)}(h,v_1)\leq \frac{s_\kappa(td(x_0,\psi(z))}{s_\kappa(d(x_0,\psi(z)))}.
\end{equation*}
Since $\md\tilde{\varphi}_{(t,z)}(1,0)=d(x_0,\psi(z))$ we finally see that
\begin{equation*}
 \jac_{k+1}(\md\tilde{\varphi}_{(t,z)})\leq d(x_0,\psi(z))\frac{[s_\kappa(td(x_0,\psi(z)))]^k}{[s_\kappa(d(x_0,\psi(z)))]^k}\jac_k(\md\psi_z)
\end{equation*}
which, together with \eqref{equation:mass-representation-param}, proves the theorem. 
\end{proof}
We finally prove that every complete ${\rm CAT}(\kappa)$-space $X$ with $\kappa<0$ admits a linear isoperimetric inequality for $\intcurr_k(X)$, $k\geq 1$.
\begin{proof}[{Proof of \thmref{Theorem:linear-isoperimetric-inequality}}]
 We fix an arbitrary $x_0\in X$ and define $S:=\cone{x_0}{T}$ as in \thmref{theorem:cat-cone-inequality}. Note that $S\in\intcurr_{k+1}(X)$ and $\bdry S=T$ by
 equation \eqref{equation:product-boundary} and the remark following it.
 For $0<r<\infty$ we have 
 \begin{equation*}
  \int_0^1\sinh^k(\sqrt{-\kappa}tr)\,dt \leq \frac{1}{kr\sqrt{-\kappa}}\sinh^k(\sqrt{-\kappa}r)
 \end{equation*}
 which, together with \thmref{theorem:cat-cone-inequality}, implies 
 \begin{equation*}
  \mass{S}\leq \frac{1}{\sqrt{-\kappa}k}\mass{T},
 \end{equation*}
 independently of the choice of $x_0$.
\end{proof}
As mentioned before the cone inequality for the case $\kappa=0$ can also be used to prove the following monotonicity formula for absolutely area minimizing 
currents, which is well-known in the case of Hadamard manifolds.
%
%
\bc\label{corollary:hadamard-monotonicity-formula}
 Let $(X,d)$ be a Hadamard space, $k\geq 1$, and let $S\in\intcurr_k(X)$ be absolutely area minimizing. If $x_0\in \spt S\ohne \spt(\bdry S)$ then the function 
 \begin{equation*}
  f:(0,\infty)\to (0,\infty)\qquad r\mapsto \frac{\|S\|(B(x_0,r))}{\omega_k r^k}
 \end{equation*}
 is monotonically non-decreasing on $(0,\dist(x_0,\spt(\bdry S))]$.
\ec
It should be mentioned that for arbitrary $S\in\intcurr_k(X)$ 
\begin{equation*}
\liminf_{r\searrow 0} \frac{\|S\|(B(x_0,r))}{\omega_k r^k}\geq 1
\end{equation*}
for $\|S\|$-almost every $x_0\in\spt S$. This follows from \cite[Theorem 9]{Kirchheim} together with the representation formula for the mass (Theorem~9.5 in 
\cite{Ambr-Kirch-curr}) and \lemref{lemma:alexandrov-metric-differential}.
\begin{proof}
 Denote by $\varrho$ the distance function to the point $x_0$ and define $\beta(r):=\|S\|(B(x_0,r))$. By the slicing theorem we have that
 $\bdry(S\rstr B(x_0,r))=\slice{S}{\varrho}{r}\in\intcurr_{k-1}(X)$ and
 \begin{equation*}
  \mass{\slice{S}{\varrho}{r}}\leq \beta'(r)
 \end{equation*}
 for almost every $r\in[0,\dist(x_0,\spt(\bdry S))]$.
 Since $S$ is absolutely area minimizing we furthermore have by \thmref{theorem:cat-cone-inequality}
 \begin{equation*}
  \beta(r)=\mass{S\rstr B(x_0,r)}\leq\mass{\cone{x_0}{\slice{S}{\varrho}{r}}}\leq\frac{r}{k+1}\beta'(r)
 \end{equation*}
 and consequently $\frac{k+1}{r}\leq \frac{d}{dt}\log(\beta(r))$ for a.e. $\,0<r<\dist(x_0,\spt(\bdry S))$. The claim now follows by integration.
\end{proof}

\section{Quasi-isometry invariance of $\divGMT_k(X)$}\label{section:quasi-isometry-invariance}
The aim of this section is to establish the quasi-isometry invariance property of the homological divergence functions stated in 
\thmref{theorem:quasi-isometry-invariance}.

Recall that two metric spaces $X$ and $Y$ are said to be quasi-isometric if there exist $\lambda\geq 1$ and $L\geq 0$ and a map $\varphi: X\to Y$ with the properties 
that
\begin{equation*}
 \frac{1}{\lambda}d(x,x') - L\leq d(\varphi(x),\varphi(x'))\leq \lambda d(x,x') + L
\end{equation*}
for all $x,x'\in X$ and that $\varphi(X)$ is $L$-dense in $Y$, i.e.~for every $y\in Y$ there exists an $x\in X$ such that $d(\varphi(x), y)\leq L$. Such a map $\varphi$ is called a quasi-isometry.

\bd
 Given two complete metric spaces $X$ and $X'$ and $\beta\in[1,\infty)$ we say that $\divGMT_k(X)\preceq_\beta \divGMT_k(X')$ if there exist $0<\varrho_0,\varrho_0'\leq 1$ and 
 $A_0,A'_0>0$ such that for every triple $(x_0, \varrho, A)$ with $x_0\in X$, $\varrho\leq\varrho_0$ and $A\geq A_0$ there exist $x_0'\in X'$, $\varrho'\leq \varrho'_0$
 and $A'\geq A'_0$ with $\delta^k_{x_0,\varrho, A}\preceq_\beta\delta^k_{x_0', \varrho', A'}$ where the two functions represent the growth functions for $X$ and $X'$, 
 respectively.
\ed 
If $\divGMT_k(X)\preceq_\beta \divGMT_k(X')$ and $\divGMT_k(X')\preceq_\beta \divGMT_k(X)$ then we write $\divGMT_k(X)\sim_\beta\divGMT_k(X')$.

In \cite[Theorem 1.1]{Brady-Farb} it was shown that $\divergence_k$ is a quasi-isometry invariant in the class of Hadamard manifolds which admit cocompact lattices.
Here, we prove the following generalization to Hadamard spaces. For this we recall that $X$ is said to be cocompact if there exists a compact set $K\subset X$ such that 
$X=\bigcup_{g\in\Gamma}gK$ where $\Gamma$ denotes the isometry group of $X$. Furthermore, $X$ is said to be proper if every closed ball of finite radius is compact.

\bt\label{theorem:generalized-quasi-isometry-invariance}
 Let $X$ and $Y$ be quasi-isometric Hadamard spaces with $Y$ proper and cocompact. Then for every $k\in\N$ we have
 \begin{equation*}
  \divGMT_k(X)\preceq_{k+1} \divGMT_k(Y).
 \end{equation*}
\et

Clearly,  \thmref{theorem:quasi-isometry-invariance} is a direct consequence of this result.
In the proof of \thmref{theorem:generalized-quasi-isometry-invariance} we will use the following special case of Theorem 1.6 in \cite{Lang-Schlichenmaier}.  
In order to state the lemma we recall that a metric space $Y$ is said to be Lipschitz 
$n$-connected for some $n\in\N$, if there exists a $c>0$ such that every $\lambda$-Lipschitz map from the boundary $\bdry\Delta^m$ of the standard $m$-dimensional simplex $\Delta^m\subset\R^{m+1}$ 
to $Y$ has a $c\lambda$-Lipschitz extension to all of $\Delta^m$ for every $m=1, \dots, n$. Note that for example Hadamard spaces and Banach spaces are Lipschitz 
$n$-connected for every $n\geq 1$.

\bl\label{lemma:Lipschitz-extension}
 Let $X$ and $Y$ be metric spaces, $r_0>0$, $\delta\geq 2$, $n\in\N$ and let $Z\subset X$ be $r_0$-separated and $\delta r_0$-dense. If the covering $\{B(z,2\delta r_0)\}_{z\in Z}$ has multiplicity at most $n$
 and if $Y$ is Lipschitz $(n-1)$-connected then every $\lambda$-Lipschitz map $\varphi: Z\to Y$ has a $C\lambda$-Lipschitz extension $\overline{\varphi}: X\to Y$ 
 of $\varphi$ for some constant $C$ depending only on $r_0$, $n$, $\delta$ and the constant from the Lipschitz connectedness.
\el

Here, we call $Z\subset X$ $r_0$-separated if $d(z,z')\geq r_0$ for all $z, z'\in Z$ with $z'\not = z$.
The lemma could be proved by applying \cite[Theorem 1.6]{Lang-Schlichenmaier}. However, since our set $Z$ is of a particularly simple form  we do not need the full strength of this
 theorem and we prefer to give a self-contained proof which follows the lines of \cite[Theorem 1.6]{Lang-Schlichenmaier} and simplifies this in our special situation.
 
\begin{proof}
 Define for each $z\in Z$ a $\frac{1}{\delta r_0}$-Lipschitz function $\nu_z: X\to \R$ by
 \begin{equation*}
  \nu_z(x):= \max\left\{0, 2 - \frac{1}{\delta r_0}d(x, z)\right\}
 \end{equation*}
 and a $\frac{2}{r_0}$-Lipschitz function $\tau_z: X\to\R$ by 
 \begin{equation*}
  \tau_z(x):= \min\left\{ \nu_z(x), \frac{2}{r_0}\dist(x, Z\ohne\{z\})\right\}.
 \end{equation*}
 Note that $\tau_z(z)=2$ and $\tau_z(z')=0$ whenever $z, z'\in Z$ and $z'\not=z$ and that for every $x\in X$ there exists a $z\in Z$ with $d(x,z)\leq \delta r_0$ and 
 $\tau_z(x)\geq 1$. Moreover, for every $x\in X$ there are at most $n$ distinct $z\in Z$ with $\tau_z(x)>0$.  Consequently, the function  
 $\overline{\tau}(x):=\sum_{z\in Z}\tau_z(x)$ is well-defined and satisfies $1\leq\overline{\tau}(x)\leq 2n$ for all $x\in X$.
 
 Define now a map $g: X\to \ell^2(Z)$ by $g(x):= (\tau_z(x)/\overline{\tau}(x))_{z\in Z}$. It can easily be checked that $g$ is $\lambda_1$-Lipschitz for a constant $\lambda_1$ only depending
 on $n$ and $r_0$. Set 
 \begin{equation*}
  \Sigma:= \left\{(v_z)\in\ell^2(Z): v_z\geq 0, \sum v_z=1\right\}
 \end{equation*}
 and denote by $\Sigma^{(m)}$ the $m$-skeleton of $\Sigma$. Observe that $g(z)=e_z\in\Sigma^{(0)}$ for every $z\in Z$, where 
 $e_z$ denotes the vertex in $\Sigma^{(0)}$ corresponding to $z\in Z$. Moreover, we have
 \begin{equation*}
  g(X)\subset\Sigma':=\left\{[e_{z_1},\dots,e_{z_k}]\subset\Sigma^{(k-1)}: k\leq n, d(z_i,z_j)<4\delta r_0\right\},
 \end{equation*}
 where $[e_{z_1},\dots,e_{z_k}]$ is the simplex spanned by these vectors.
 For $z\in Z$ define $h^{(0)}(e_z):= \varphi(z)$ and extend $h^{(0)}$ recursively to  maps $h^{(m)}:\Sigma'\cap\Sigma^{(m)}\to Y$ for $m=1,\dots, n-1$ in the obvious way using the Lipschitz connectedness
 of $Y$. Clearly, $h:= h^{(n)}$ is $\lambda_2$-Lipschitz on every closed simplex of $\Sigma'$ for some $\lambda_2$ only depending on $\lambda$, $r_0$, $\delta$, $n$ and the  constant from the Lipschitz connectedness.
 Finally, set $\overline{\varphi}:= h\circ g$ and note that this is an extension of $\varphi$ since $\overline{\varphi}(z)= h(g(z))=h(e_z)=\varphi(z)$.
 
 We are left to check that $\overline{\varphi}$ is Lipschitz.
 For this let $x,x'\in X$ with $x'\not= x$. Let $S, S'\subset \Sigma'$ be the minimal simplices containing $g(x)$ and $g(x')$, respectively.
 If $d(x,x')\geq \frac{1}{4}r_0$ then choose $z, z'\in Z$ with $d(x,z), d(x',z')\leq \delta r_0$ and such that $\tau_z(x)\geq 1$ and $\tau_{z'}(x')\geq 1$ and note that 
 $g(z)\in S$ and $g(z')\in S'$. We now simply calculate
 \begin{equation*}
  \begin{split}
   d(\overline{\varphi}(x), \overline{\varphi}(x')) &\leq d(\overline{\varphi}(x), \overline{\varphi}(z)) + d(\overline{\varphi}(z), \overline{\varphi}(z')) + d(\overline{\varphi}(z'), \overline{\varphi}(x'))\\
    &\leq \lambda_1\lambda_2 d(x,z) +\lambda d(z, z') + \lambda_1\lambda_2 d(x', z')\\
    &\leq (\lambda_1\lambda_2+\lambda)[d(x,z)+d(x',z')] +\lambda d(x,x')\\
    &\leq \left[8(\lambda_1\lambda_2 + \lambda)\delta+\lambda\right]d(x,x').
  \end{split}
 \end{equation*}
 On the other hand, if $d(x,x')<\frac{1}{4}r_0$ then choose $z\in Z$ with $d(x,z)\leq\delta r_0$ and $\tau_z(x)\geq 1$ and observe that this implies $\tau_z(x')>0$. In particular, $S\cap S'\not=\emptyset$ and therefore there 
 exists a $v\in S\cap S'$ such that
 \begin{equation*}
  d(g(x), v)+d(v, g(x'))\leq c'd(g(x), g(x'))
 \end{equation*}
 for a constant $c'$ depending only on $n$ and hence
 \begin{equation*}
  \begin{split}
   d(\overline{\varphi}(x),\overline{\varphi}(x')) &\leq  d(\overline{\varphi}(x), h(v)) + d(h(v), \overline{\varphi}(x'))\\
    &\leq \lambda_2 d(g(x), v) + \lambda_2 d(v,g(x'))\\
    &\leq \lambda_1\lambda_2 c' d(x,x').
  \end{split}
 \end{equation*}
 This concludes the proof.
\end{proof}

We note the following simple but important instance where the hypotheses are satisfied.

\bc
 Let $X$ be a proper and cocompact metric space, $Y$ a Hadamard space and $\varphi: X\to Y$ an $(\lambda, L)$-quasi-isometry for some $\lambda\geq 1$ and $L>0$. 
 Then there exists a $2\lambda L$-separated
 and $4\lambda L$-dense set $Z\subset X$ and a Lipschitz extension $\overline{\varphi}:X\to Y$ of the $2\lambda$-biLipschitz map $\varphi\on{Z}$.
\ec

\begin{proof}
 Set $r_0:= 2\lambda L$ and note that by the hypotheses on $X$ there exists an $n\in\N$ such that every closed ball of radius $4r_0$ in $X$ contains at most $n$
 pairwise disjoint open balls of radius $r_0/2$. Choose a maximally $r_0$-separated set $Z\subset X$. Then $Z$ is $2r_0$-dense and the covering 
 $\{B(z,4r_0)\}_{z\in Z}$ has multiplicity at most $n$. Since $Y$ is Lipschitz $(n-1)$-connected the claim now follows directly from \lemref{lemma:Lipschitz-extension}.
\end{proof}

Finally, we are in a position to prove the quasi-isometry invariance of $\divGMT_k$.

\begin{proof}[{Proof of \thmref{theorem:generalized-quasi-isometry-invariance}}]
  Let $X$ and $Y$ be as in the hypotheses and let $\varphi: X\to Y$ be a $(\lambda, L)$-quasi isometry for some $\lambda\geq 1$ and $L\geq 0$. We may assume without loss of
  generality that $L>0$. Set $r_0:= 2\lambda L$.
  Set moreover $\varrho':=\varrho'_0:=\frac{1}{2}$, choose $\varrho_0\in(0,1)$ small enough (as below) and let $A_0, A'_0\geq 1$ be arbitrary. Choose
  $x_0\in X$, $0<\varrho\leq \varrho_0$, and $A\geq A_0$ arbitrarily.
  Let $r>0$ be large enough, to be specified below, and choose $T\in \intcurr_k(X)$ with  $\spt T\subset S(x_0, r)$ compact and such that $\bdry T=0$ and 
  $\mass{T}\leq Ar^k$. 
  By \cite[Theorem 1.6]{Wenger-GAFA} there exists an $R\in\intcurr_{k+1}(X)$ which is absolutely area minimizing and satisfies $\bdry R=T$. 
  Set  $\Omega:=\spt R\cap B(x_0, \frac{7}{8}r)$, note that $\Omega$ is compact, and let $Z\subset\Omega$ be maximally $r_0$-separated. 
  Extend $Z$ to a  
  maximally $r_0$-separated subset $\hat{Z}$ of $X$ and note that $\varphi\on{\hat{Z}}$ is $2\lambda$-biLipschitz and that $Z':=\varphi(\hat{Z})$ is 
  $r'_0$-separated and $\delta' r'_0$-dense in $Y$, where $r'_0:= L$ and $\delta':= 2\lambda+1$.
  Next, we verify that for all $x\in\Omega$ and for all $s\in(0,\frac{1}{8}r)$
  \begin{equation}\label{equation:Ahlfors-support}
   D_1s^{k+1}\leq \|R\|(B(x,s))\leq D_2s^{k+1}
  \end{equation}
  for constants $0< D_1\leq D_2<\infty$ depending only on the isoperimetric constant $C$ for $\intcurr_k(X)$ and on $A$ and $k$. Indeed, the first inequality 
  follows directly from 
  \lemref{Lemma:Q-min-growth} and \remref{remark:absolute-minimizers-growth-also}
  whereas for the second inequality it is enough to note that by \corref{corollary:hadamard-monotonicity-formula} and by the 
  isoperimetric inequality for $\intcurr_k(X)$ we have
  \begin{equation*}
   \frac{\|R\|(B(x, s))}{s^{k+1}}\leq \frac{\|R\|(B(x, \frac{1}{8}r))}{(\frac{1}{8}r)^{k+1}}\leq \left(\frac{8}{r}\right)^{k+1}\mass{R}\leq 8^{k+1}CA^{\frac{k+1}{k}}.
  \end{equation*}
 As a direct consequence of \eqref{equation:Ahlfors-support} we obtain that the covering $\{B(z,4r_0)\cap\Omega\}_{z\in Z}$ of $\Omega$ has multiplicity at most 
 $n:= \frac{9^{k+1}D_2}{D_1}$. 
 By \lemref{lemma:Lipschitz-extension} there thus exists a
 $\lambda_1$-Lipschitz extension $\overline{\varphi}: \Omega\to Y$ of $\varphi\on{Z}$ for some constant $\lambda_1$ only depending on $\lambda$, $r_0$ and $n$.
 
 Now, it follows as in the corollary above that the covering $\{B(z', 2\delta' r'_0)\}_{z'\in Z'}$ of $Y$ has multiplicity at most $n'$ for some finite number $n'\in\N$.
 Thus, by \lemref{lemma:Lipschitz-extension}, there exists a
 $\lambda_2$-Lipschitz extension $\overline{\eta}: Y\to X$ of $\eta:= (\varphi\on{\hat{Z}})^{-1}: Z'\to X$ for some constant $\lambda_2$ only depending on $\lambda$, $\delta'$,
 $r'_0$ and $n'$.
 It follows that 
 \begin{equation}\label{equation:distance-geodesics-projection}
  d(x,(\overline{\eta}\circ\overline{\varphi})(x))\leq a_1:=2(1+\lambda_1\lambda_2)r_0
 \end{equation}
 for all $x\in\Omega$ and that 
 \begin{equation}\label{equation:back-distortion}
  d(x'_0, y)\leq \lambda d(x_0,\overline{\eta}(y)) + a_2
 \end{equation}
 for all $y\in Y$ where $x'_0:=\varphi(x_0)$ and $a_2:= (1+\lambda\lambda_2)\delta'r'_0+L$.
 Let $\tau:[0,1]\times X\to X$ be the Lipschitz map for which $t\mapsto\tau(t, x)$ is the geodesic from $x$ to 
 $(\overline{\eta}\circ\overline{\varphi})(x)$ parametrized proportional to arc-length.
 Define $\pi:\Omega\to \R$ by $\pi(x):= d(x'_0, \overline{\varphi}(x))$ and observe that $\pi$ is $\lambda_1$-Lipschitz and satisfies
 \begin{equation}\label{equation:distortion-qi}
  \frac{1}{\lambda}d(x,x_0) - L_1\leq \pi(x)\leq \lambda d(x,x_0) + L_1
 \end{equation}
 for all $x\in\Omega$ where $L_1:=2(\lambda+\lambda_1)r_0+L$.
 Set furthermore $r_1:=\frac{6}{8\lambda}r-L_1$ and $r_2:=\frac{7}{8\lambda}r-L_1$ and note that for $\Omega_t:= \{x\in \Omega: \pi(x)\leq t\}$ we have
 \begin{equation}\label{equation:subset-inclusions-Omega}
 \Omega_t\subset B\left(x_0, \lambda (t+L_1)\right)\quad\text{ and }\quad \Omega\ohne\Omega_t\subset X\ohne B\left(x_0, \frac{1}{\lambda}(t-L_1)\right)
 \end{equation}
 for all $t\geq L_1$, by \eqref{equation:distortion-qi}. 
 By the slicing theorem there exists an $r'\in (r_1, r_2)$ such that $R\rstr \Omega_{r'}\in \intcurr_{k+1}(X)$,
 \begin{equation*}
  \spt \bdry (R\rstr\Omega_{r'})\subset\{x\in\Omega: \pi(x) = r'\}\quad\text{ and }\quad
  \mass{\bdry(R\rstr \Omega_{r'})}\leq \tilde{A}r^k,
 \end{equation*}
 where $\tilde{A}=8\lambda\lambda_1 CA^{\frac{k+1}{k}}$. Indeed, the first and second claim follow directly from 
 the slicing theorem whereas for the third claim it is enough to note that again by the slicing theorem
 $R\rstr \Omega_t\in\intcurr_{k+1}(X)$ for almost every  $t\in(r_1, r_2)$ and 
 \begin{equation*}
  \int_{r_1}^{r_2}\mass{\bdry(R\rstr \Omega_t)}dt\leq \lambda_1\int_{r_1}^{r_2}\frac{d}{dt}\|R\|(\Omega_t)dt\leq \lambda_1\mass{R}\leq \lambda_1CA^{\frac{k+1}{k}}r^{k+1}.
 \end{equation*}
 Set $T':= \overline{\varphi}_\# (\bdry (R\rstr \Omega_{r'}))$ and observe that $T'\in\intcurr_k(Y)$ and $\bdry T'=0$. Furthermore, we have
 $\spt T'\subset S(x'_0, r')$ and $\mass{T'}\leq A'r^k$, where $A':=\lambda_2^k\tilde{A}$.
 
 Now choose $S'\in\intcurr_{k+1}(Y)$ such that $\bdry S'=T'$ and $\spt S'\subset Y\ohne U(x'_0, \varrho'r')$ as well as 
 $\mass{S'}\leq 2\delta^k_{x'_0, \varrho', A'}(r').$
 Set
 \begin{equation*}
  S:= R- R\rstr \Omega_{r'} - \tau_\#([0,1]\times \bdry (R\rstr \Omega_{r'})) + \overline{\eta}_\# S'.
 \end{equation*}
 It is clear that $S\in\intcurr_{k+1}(X)$. Furthermore, we compute
 \begin{equation*}
  \bdry S = T - \bdry(R\rstr \Omega_{r'}) - \tau_\#([1]\times \bdry (R\rstr \Omega_{r'}) - [0]\times \bdry (R\rstr B_{r'})) + \overline{\eta}_\# T'
          = T
 \end{equation*}
 and
 \begin{equation*}
  \begin{split}
   \mass{S}&\leq \lambda_2^{k+1}\mass{S'} + \mass{R} + a_1(k+1)(\lambda_1\lambda_2)^k\mass{\bdry (R\rstr \Omega_{r'})}\\
           &\leq 2\lambda_2^{k+1}\delta^k_{x'_0,\varrho', A'}(r') + CA^{\frac{k+1}{k}}r^{k+1} + a_1(k+1)(\lambda_1\lambda_2)^k\tilde{A}r^k.
  \end{split}
 \end{equation*}
 Finally, we have $\spt S\subset X\ohne U(x_0, \varrho r)$ since by \eqref{equation:distance-geodesics-projection} and  \eqref{equation:subset-inclusions-Omega} 
 \begin{equation*}
  \spt(R - R\rstr\Omega_{r'}) \subset X\ohne U\left(x_0, \frac{1}{\lambda}(r'-L_1)\right)
 \end{equation*}
 and
 \begin{equation*}
  \spt(\tau_\#([0,1]\times\bdry(R\rstr\Omega_{r'}))) \subset X\ohne U\left(x_0, \frac{1}{\lambda}(r'-L_1)-a_1\right)
 \end{equation*}
 and by \eqref{equation:back-distortion}
 \begin{equation*}
  \spt\left(\overline{\eta}_\# S'\right)\subset X\ohne U\left(x_0, \frac{1}{\lambda}(\varrho'r' - a_2)\right).
 \end{equation*}
 This concludes the proof provided that $r$ was chosen large enough and $\varrho_0$ small enough.
\end{proof}

\end{document}